\documentclass[11 pt]{article}  
\usepackage[utf8]{inputenc}
\usepackage{amsmath}
\usepackage{amsfonts}
\usepackage{amssymb}
\usepackage{graphicx}
\usepackage{mathrsfs}
\usepackage{upref,amsthm,amsxtra,exscale}
\usepackage{stmaryrd}
\usepackage{cite}
\usepackage[colorlinks=true,urlcolor=blue,
citecolor=red,linkcolor=blue,linktocpage,pdfpagelabels,
bookmarksnumbered,bookmarksopen]{hyperref}

\usepackage{fullpage}

\usepackage{subcaption}
\usepackage{caption}
\usepackage{cleveref}

\usepackage{enumitem}

\newtheorem{theorem}{Theorem}[section]

\newtheorem{lemma}[theorem]{Lemma}
\newtheorem{proposition}[theorem]{Proposition}

\numberwithin{equation}{section}

\theoremstyle{definition}
\newtheorem{remark}[theorem]{Remark}

\usepackage{color}
\usepackage[dvipsnames]{xcolor}

\usepackage[T1]{fontenc}
\usepackage[utf8]{inputenc}
\usepackage{authblk}

\title{\textbf{Existence and Multiplicity of Normalized Solutions to Schr\"{o}dinger Equations with General Nonlinearities in Bounded Domains}}
\author[a,b]{Wei Ji\thanks{Email: jiwei2020@amss.ac.cn} }
\affil[a]{\small Academy of Mathematics and Systems Science, Chinese Academy of Sciences, Beijing 100190, PR China}
\affil[b]{\small{University of Chinese Academy of Sciences, Beijing 100049, P.R. China}}
  \date{}

\begin{document}
  
  \maketitle

\abstract{This paper focuses on the existence of multiple normalized solutions 
to Schr\"{o}dinger equations with general nonlinearities in bounded domains via variational methods. 
We first obtain two positive normalized solutions, one is a normalized ground state by searching for a local minimizer, and the other one is a mountain pass solution. 
Secondly, using a version of Linking theorems for normalized solutions, we prove the multiplicity of solutions to Schr\"{o}dinger equations in a star-shaped bounded domain.
Moreover, we arrive at the existence of nonradial normalized solutions to Schr\"odinger equations in a ball.  
 }

\bigskip
\bigskip

\noindent\text{\textbf{Keywords:}}  Normalized solutions, positive solutions,
 multiple solutions, nonradial solutions, general nonlinearities.
\medskip

\noindent\text{\textbf{MSC2020:}} 35A01, 35B09, 35A15, 35B40, 35J60.

\medskip

\section{Introduction and main results}

In this paper, we study the existence of positive and multiple normalized solutions for the semi-linear Dirichlet problem:
\begin{align}
\label{Peps}
\begin{cases}-\Delta u +\lambda u=f(u) \quad\mbox{ in }\Omega,\qquad
  \\ u|_{\partial \Omega} =0\end{cases}
\end{align}
with prescribed $L^2$-norm
\begin{equation}
  \int_{\Omega}\vert u \vert^2 =c,
\end{equation}
where $\Omega$ is a bounded, smooth and star-shaped domain in $\mathbb R^N$ and $N \geq 3$.

The normalized solutions to nonlinear Schr\"{o}dinger equation~\eqref{Peps} have been of constant attention for many years.
In the case $\Omega$ is a ball and $f(u)=|u|^{p-2}u$ in ~\eqref{Peps} where $2<p<2^*$, 
the authors in~\cite{noris2015existence} focus a two-constraint problem, i.e.
$$
\max\left\{\int_{\Omega}|u|^pdx:u\in H_0^1(\Omega), \int_{\Omega}u^2dx=1, \int_{\Omega}|\nabla u|^2dx=a\right\}, 
$$
to establish a global branch respect to $\lambda$ of positive solution of ~\eqref{Peps} relying on the uniqueness results in~\cite{zhang1992uniqueness}, 
and then obtain the existence and nonexistence of positive normalized solutions. 
Recently, the authors in \cite{songzou2024two} study two positive normalized solution by searching for a local minimizer and a mountain pass solution
for Brezis-Nirenberg problem. And the authors in~\cite{Qishijie} obtain the multiple normalized solutions of\eqref{Peps} when the nonlinearity 
is sobolev subcritical and nonhomogeneous, by establishing special links and using the deformation method on the mass constraint manifold.
 Furthermore, they actually consider
the non-autonomous equation with potentials.

In this paper, first, we consider the existence of positive normalized solutions with mass supercritical general nonlinearities. 
We assume on $f$:
\begin{itemize}
  \item [$(f_1)$]   $f(t)\in C( \mathbb R, \mathbb R)$, $f(0)=0$, and there exist $\alpha,\beta$ satisfying $2+\frac{4}{N}<\alpha\leq\beta<2^*$, where $2^*=\frac{2N}{N-2}$, 
  such that
  \begin{align}\label{estimateB}
    0<\alpha F(u)\leq f(u)u\leq\beta F(u)\quad, u\neq0, 
    \end{align}
    where $F(u)=\int_0^uf(u)ds.$
  \end{itemize}
By ~\eqref{estimateB}, we can deduce that there exist $\mu, \zeta$ satisfies $0<\mu<\zeta$ such that
\begin{align}\label{estimateC}
\mu (|t|^{ \alpha } +|t|^{\beta })\leq F(t)\leq \zeta (|t|^{ \alpha} +|t|^{\beta }). 
\end{align}

To find the positive normalized solutions of \eqref{Peps}, we search for critical points of the energy
\begin{align}
  E(u)=\frac12\int_\Omega|\nabla u|^2dx-\int_{\Omega}F_+(u)dx,
\end{align}

under the constraint
$$\int_\Omega|u^+|^2dx=c,$$
where $u^+=\max\{u,0\}$,  and $F_+(t)=\int_{0}^{t}f_+(s)ds$, where
$f_+(t)$ is defined by
\begin{align}
  f_+(t)=
\begin{cases}
  f(u), \quad  t \geq 0,
    \\ 0, \quad \quad ~~ t<0. 
  \end{cases}
  \end{align}

Indeed, by ~\eqref{estimateB}, we know that $f(t)>0$ for $t>0$, $f(t)<0$ for $t<0$ and $F(t)>0$ for $t \neq 0$. 
Therefore, $f_+(u)=f(u^+)$ and $F(u) \geq F_+(u)=F(u^+)$.

We set
$$
S_c^+:=\left \{ u\in H_0^1(\Omega):\int_\Omega|u^+|^2dx=c \right\}.
$$

By~\cite{willem1996minimax}, any critical point $u$ of $E|_{S_c^+}$ satisfies the following Pohozaev identity:
$$
\int_{\Omega}|\nabla u|^2dx-\frac{1}{2}\int_{\partial\Omega}|\nabla u|^2\sigma\cdot nd\sigma=
\frac{N}{2}\int_{\Omega}f_+(u)u^+dx-N\int_{\Omega}F_+(u)dx.
$$

Note that $\sigma\cdot n>0$ since $\Omega$ is star-shaped with respect to the origin. Hence, $u$ belongs

to $\mathcal{G}$ where
$$
\mathcal{G}:=\left\{u\in S_c^+:\int_\Omega|\nabla u|^2dx>\frac{N}{2}\int_{\Omega}f_+(u)u^+dx-N\int_{\Omega}F_+(u)dx\right\}.
$$

We will prove that $\mathcal{G}$ is nonempty and the lower bound of $ E(u)$ on $\mathcal{G}$ can be obtained. As a consequence,
we get a normalized ground state of \eqref{Peps}. Furthermore, by establishing a mountain pass structure, we obtain another positive solution.

\medskip

Our main conclusions are as follows. 
 
\begin{theorem}\label{theorem1}
Let $N \geq 3$, $\Omega$ be bounded, smooth and star-shaped domain with respect to the origin and $f$ satisfy $(f_1)$.  
Then for any 
\begin{align}\label{upboundforc}
  c<\sup_{u\in S_c^+}\left(min\{c(g_{1,u}), (\frac{1}{2}-\frac{2}{N(\alpha-2)})C_1g_{2, u}^{-1}\}\right)
\end{align}

where $c(g_{1,u})$ satisfies $\frac{(\beta-2)\zeta N}{2}(c(g_{1,u})^{\frac{\alpha-2}{2}}+c(g_{1,u})^{\frac{\beta-2}{2}})=g_{1,u}$, 
$$
g_{1,u}=\frac{\int_{\Omega}|\nabla u|^2dx}{\int_{\Omega}|u^+|^{\alpha}+|u^+|^{\beta}dx}, ~~g_{2,u}=\frac{1}{2}\int_{\Omega}|\nabla u|^2dx,
$$
and $C_1$ is a constant defined in~\eqref{equaofc1}, 
equation \eqref{Peps} admits a positive normalized solution $(\lambda_c, u_c)$ such that $\int_{\Omega}|u_c|^2=c$. Moreover, $u_c$ is 
a normalized ground state of \eqref{Peps}. 

\end{theorem}

\begin{theorem}\label{theorem2}
  Let $N \geq 3$, $\Omega$ be bounded, smooth and star-shaped domain with respect to the origin, $f$ satisfy $(f_1)$ and 
  $c$ satisfy~\eqref{upboundforc}. Then equation \eqref{Peps} admits a normalized solution $(\widetilde{\lambda}_c, \widetilde{u}_c)$ 
  such that $\int_{\Omega}|\widetilde{u}_c|^2=c$ and $\widetilde{u}_c \neq u_c$. 
\end{theorem}

Furthermore, we supplement some results of positive normalized solutions of~\eqref{Peps} when $f$ is combined with mass supercritical and
critical or subcritical terms. 

\begin{theorem}\label{theorem3}
  Let $N \geq 3$, $\Omega$ be bounded, smooth and star-shaped domain with respect to the origin, $f=|u|^{p-2}u+a|u|^{q-2}u$, 
  where $2<q\leq 2+\frac{4}{N}<p<2^*$ and 
  $c$ satisfies
 
\begin{align}\label{upuqccage0}
c<\sup_{u\in S_c^+}\left\{min\{c_1, c(g_{3,u}),(\frac{1}{4}-\frac{1}{N(p-2)})C_2g_{2,u}^{-1}\}\right\}    
\end{align}

if $a>0$, 
where $c_1$ satisfies
$$
\frac{a(p-q)}{q(p-2)}C_{N,q}^q c_1^{\frac{2q-N(q-2)}4}
=(\frac{1}{4}-\frac{1}{N(p-2)})(C_2)^{\frac{4-N(q-2)}4},
$$
$c(g_{3,u})$ satisfies
$$
(\frac{1}{2}-\frac{1}{p})c(g_{3,u})^{\frac{p-2}{2}}+a(\frac{1}{2}-\frac{1}{q})c(g_{3,u})^{\frac{q-2}{2}}=g_{3,u},
$$ 
where 
$$
g_{3,u}=\frac{\int_{\Omega}|\nabla u|^2dx}{\int_{\Omega}|u^+|^{p}dx+a\int_{\Omega}|u^+|^{q}dx},
$$
and $C_2$ is defined in ~\eqref{defineC2}; or
\begin{align}\label{upuqccale0}
c<\sup_{u\in S_c^+}\left\{min\{c(g_{4,u}), c_2(g_{2,u},g_{5,u})\}\right\}
\end{align}
if $a<0$, where
 $c(g_{4,u})$ satisfies $(\frac{1}{2}-\frac{1}{p})c(g_{4,u})^{\frac{p-2}{2}}=g_{4,u}$, 
 where
$$
g_{4,u}=\frac{\int_{\Omega}|\nabla u|^2dx}{\int_{\Omega}|u^+|^{p}dx}, 
$$
$$
c_2(g_2(u),g_5(u))= min\left\{\varsigma_0(\frac{1}{2}-\frac{2}{N(p-2)})C_3g_{2,u}^{-1},~
\left ((1-\varsigma_0)\frac{q}{|a|}(\frac{1}{2}-\frac{2}{N(p-2)})C_3g_{5,u}\right )^{\frac{2}{q}}\right\},
$$
where
$$
g_{5,u}=\frac{1}{\int_{\Omega}|u^+|^{q}dx}, 
$$
$\varsigma_0$ is defined in \eqref{definevarsigma0}
and $C_3$ is defined in ~\eqref{defineC3}.

  Then \eqref{Peps} admits two normalized solution $(\lambda_c, u_c)$ and
  $(\widetilde{\lambda}_c, \widetilde{u}_c)$ 
  such that $\widetilde{u}_c \neq u_c$ and $\int_{\Omega}|u_c|^2=\int_{\Omega}|\widetilde{u}_c|^2=c$. 
  Moreover, $u_c$ is a normalized ground state and $\widetilde{u}_c$ is a mountain pass solution. 
\end{theorem}

\begin{remark}
In Theorem~\ref{theorem1}, since $\frac{\alpha}{2}, \frac{\beta}{2}>1$ and $\int_{\Omega} |\nabla u|^2 >\lambda_1 c>0$, 
where $\lambda_1$ denotes the first Dirichlet eigenvalue of $-\Delta$ on $\Omega$ and it is standard to know that $\lambda_1>0$.
then $c(g_{1,u})$ is well defined, and we can verify that
the range of $c$ is suitable. 
Similar results apply to Theorem~\ref{theorem2} and Theorem~\ref{theorem3} as well. 
\end{remark}

\medskip

To find the general normalized solutions of \eqref{Peps}, we search for critical points of the energy
  \begin{align}\label{energyJ}
   \widetilde{E}(u) = \frac12\int_\Omega|\nabla u|^2dx-\int_\Omega F(u).
  \end{align}
We have the following results. 
\begin{theorem}\label{theorem4A}
  Let $N \geq 3$, $\Omega$ be bounded, smooth and star-shaped domain with respect to the origin, $f$ satisfy $(f_1)$, 
   and $c<c_k$ where $c_k$ is defined by~\eqref{fdefineck}. 
  Then equation \eqref{Peps} admits $k$ normalized solutions.  
  \end{theorem}

\begin{theorem}\label{theorem4B}
Let $N \geq 3$, $\Omega$ be bounded, smooth and star-shaped domain with respect to the origin, $f=|u|^{p-2}u-|u|^{q-2}u$, 
where $2<q<p$ and $2+\frac{4}{N}<p<2^*$, and $c<\alpha_k$ where $\alpha_k$ is defined by~\eqref{up-uqdefineck}.  
Then equation \eqref{Peps} admits $k$ normalized solutions.  
\end{theorem}

\begin{theorem}\label{theorem5}
  Let $N \geq 4$, $\Omega=B$ be a ball, $f$ satisfy $(f_1)$, and $c<c_k$ where $c_k$ is defined by~\eqref{fdefineck}(when $\Omega=B$),  
  Then equation \eqref{Peps} admits $k$ nonradial sign-changing normalized solutions.  
  \end{theorem}

\medskip

The rest of this paper is organized as follows: In section 2, we focus on the normalized ground state, that is, proving Theorem~\ref{theorem1}.
Section 3 is devoted to the mountain pass solution. In section 4, we supplement some results of mixed nonliearities. Finally, we finish this paper by studying the nonradial sign-changing normalized solutions in a ball.

  \section{A normalized ground state}

In this section, we study the local minimizer of $E(u)$ on $\mathcal{G}$ and obtain a
positive normalized solution of~\eqref{Peps}, furthermore, proved to be a ground state of~\eqref{Peps}.

\begin{lemma}\label{EinGandpartialG}
Let $N \geq 3$, $\Omega$ be bounded, smooth and star-shaped domain with respect to the origin and $f$ satisfy $(f_1)$. Then we have 
  $\inf_{u \in \mathcal{G}}E(u)>0$, and
   any sequence $\{u_n\} \subset \mathcal{G}$ satisfying $\limsup_{n\to \infty} E(u_n) < +\infty$ is bounded in $H_0^1(\Omega)$.

  \end{lemma} 

\begin{proof}
For any $u\in \mathcal{G}$, by~\eqref{estimateB}, we have
    $$
    \frac{( \alpha -2)N}{2}\int_{\Omega}F_+(u)dx \leq \frac{N}{2}\int_{\Omega}f_+(u)u^+dx-N\int_{\Omega}F_+(u)dx
    \leq\frac{(\beta -2)N}{2}\int_{\Omega}F_+(u)dx,
    $$
then
$$
 \begin{aligned}
  \frac{2}{(\beta -2)}(\frac{1}{2}\int_{\Omega}f_+(u)u^+dx-\int_{\Omega}F_+(u))dx &\leq \int_{\Omega}F_+(u)dx\\
  &\leq \frac{2}{( \alpha -2)}(\frac{1}{2}\int_{\Omega}f_+(u)u^+dx-\int_{\Omega}F_+(u))dx.
\end{aligned}
$$

    Therefore, 
    \begin{align}\label{eugradientu}
    E(u)=\frac12\int_\Omega|\nabla u|^2dx-\int_{\Omega}F_+(u)dx
    >(\frac12-\frac{2}{( \alpha -2)N})\int_\Omega|\nabla u|^2dx.
    \end{align}
   
    Note that 
    \begin{align}\label{nablaugeqlambda1c}
     \int_\Omega|\nabla u|^2dx \geq \lambda_1c>0. 
    \end{align}
   
   From~\eqref{eugradientu} and~\eqref{nablaugeqlambda1c} we derive that $\inf_{u \in \mathcal{G}} E(u)>0$ on $\mathcal{G}$.  
   
    Let $\{u_n\} \subset \mathcal{G}$ satisfying $\limsup_{n\to \infty} E(u_n) < +\infty$. 
   We have
   $$
   \limsup_{n \to \infty}(\frac12-\frac{2}{( \alpha -2)N})\int_\Omega|\nabla u|^2dx<\limsup_{n \to \infty}E(u_n)<+\infty. 
   $$
   Therefore, $\{u_n\} \subset \mathcal{G}$ is bounded in $H_0^1(\Omega)$. 
\end{proof}

\begin{lemma}\label{EinGandpartialG2}
Under the hypotheses of Lemma~\ref{EinGandpartialG}, assume that 
\eqref{upboundforc} holds true, then 
 $\mathcal{G} \neq \varnothing$, and 
$$
0<\inf_{u\in \mathcal{G}} E(u)<\inf_{u\in \partial \mathcal{G}}E(u).
$$

 \end{lemma}

 \begin{proof}
First, for any $2<p<2^*$, by Sobolev inequality, we have
\begin{align}\label{estimateforup}
\int_\Omega|u|^p \leq (S_{p^*}^{-1}\int_\Omega|\nabla u|^{p^*}dx)^{\frac{p}{p^*}}, 
\end{align}

where $S_{p^*}$ is the Sobolev optimal constant with respect to $p^*$ and $p^*=\frac{Np}{N+p}$. Since $2<p<2^*$, then $1<\frac{2N}{N+2}<p^*<2$.

Furthermore, by H\"older inequality, we have
\begin{align}\label{estimatefornablau}
\int_\Omega|\nabla u|^{p^*}dx \leq (\int_\Omega(|\nabla u|^{p^*})^{\frac{2}{p^*}}dx)^{\frac{p^*}{2}}\cdot |\Omega|^{\frac{2-p^*}{2}}. 
\end{align}
 
Combining~\eqref{estimateforup} and~\eqref{estimatefornablau}, we have
\begin{align}\label{sobinequality}
  \int_\Omega|u|^p \leq (S_{p^*}^{-1}|\Omega|^{\frac{2-p^*}{2}})^{\frac{p}{p^*}} (\int_\Omega|\nabla u|^2dx)^{\frac{p}{2}}.
\end{align}

Since any $u \in \partial\mathcal{G}$ satisfies
$$
\int_\Omega|\nabla u|^2dx=\frac{N}{2}\int_{\Omega}f_+(u)u^+dx-N\int_{\Omega}F_+(u)dx, 
$$

then

$$
\begin{aligned}
\int_\Omega|\nabla u|^{2}dx &\leq \frac{(\beta-2)N}{2}\int_\Omega F(u)dx \\
  &\leq \frac{(\beta-2)N\zeta }{2} (\int_\Omega |u|^{\alpha}dx+\int_\Omega|u|^{\beta}dx)\\
  &\leq \frac{(\beta-2)N\zeta }{2}\left((S_{\alpha^*}^{-1}|\Omega|^{\frac{2-\alpha^*}{2}})^{\frac{\alpha}{\alpha^*}}
   (\int_\Omega|\nabla u|^2dx)^{\frac{\alpha}{2}}
   +(S_{\beta^*}^{-1}|\Omega|^{\frac{2-\beta^*}{2}})^{\frac{\beta}{\beta^*}} (\int_\Omega|\nabla u|^2dx)^{\frac{\beta}{2}}\right). 
  \end{aligned}
$$

We know that $\frac{\alpha}{2}, \frac{\beta}{2}>1$, 
then there exist $C_1=C_1(N, \alpha,\beta)>0$ such that
\begin{align}\label{GradientuandCp}
  \int_\Omega|\nabla u|^2dx \geq C_1. 
\end{align}
Indeed, $C_1$ satisfies
\begin{align}\label{equaofc1}
    \frac{(\beta-2)N\zeta }{2}\left((S_{\alpha^*}^{-1}|\Omega|^{\frac{2-\alpha^*}{2}})^{\frac{\alpha}{\alpha^*}}
    (C_1)^{\frac{\alpha}{2}-1}
    +(S_{\beta^*}^{-1}|\Omega|^{\frac{2-\beta^*}{2}})^{\frac{\beta}{\beta^*}} (C_1)^{\frac{\beta}{2}-1}\right)=1. 
\end{align}
It is not difficult to verify that $C_1>0$ satisfying~\eqref{equaofc1} is unique.

Therefore, 
$$
\inf_{u \in \partial\mathcal{G}}E(u)\geq (\frac{1}{2}-\frac{2}{N(\alpha-2)})\int_\Omega|\nabla u|^2dx
\geq (\frac{1}{2}-\frac{2}{N(\alpha-2)})C_1.
$$

Secondly, for any $u\in S_1^+$ and $c$ satisfying $\frac{(\beta-2)\zeta N}{2}(c^{\frac{\alpha-2}{2}}+c^{\frac{\beta-2}{2}})<g_{1,u}$, where
$$
g_{1,u}=\frac{\int_{\Omega}|\nabla u|^2dx}{\int_{\Omega}|u^+|^{\alpha}+|u^+|^{\beta}dx},
$$

$$
\begin{aligned}
    \int_\Omega|\nabla (\sqrt{c}u)|^{2}dx &=c\int_\Omega|\nabla u|^{2}dx \\
  &> \frac{(\beta-2)\zeta N}{2}(c^{\frac{\alpha}{2}}\int_\Omega| u^+dx|^{\alpha}+c^{\frac{\beta}{2}}\int_\Omega| u^+|^{\beta}dx)\\
  &=\frac{(\beta-2)\zeta N}{2}(\int_\Omega| cu^+|^{\alpha}dx+\int_\Omega| cu^+|^{\beta}dx)\\
  &\geq \frac{(\beta-2) N}{2}\int_\Omega|F_+(\sqrt{c}u)dx\\
  &\geq \frac{N}{2}\int_{\Omega}f_+(\sqrt{c}u)\sqrt{c}u^+dx-N\int_{\Omega}F_+(\sqrt{c}u)dx
  \end{aligned}
$$

Hence $\sqrt{c}u \in \mathcal{G}$ and $\mathcal{G}$ is not empty.

Lastly, for any $u\in S_1^+$ and $c<(\frac{1}{2}-\frac{2}{N(\alpha-2)})C_1g_{2, u}^{-1}$, where
 $$
 g_{2,u}=\frac{1}{2}\int_{\Omega}|\nabla u|^2 dx, 
$$
 we have
$$
E(\sqrt{c}u)=\frac{c}{2}\int_\Omega|\nabla u|^{2}dx-\int_\Omega|F^+(\sqrt{c}u^+)
<\frac{c}{2}\int_\Omega|\nabla u|^{2}dx<(\frac{1}{2}-\frac{2}{N(\alpha-2)})C_1.
$$

When \eqref{upboundforc} holds, we can take $u\in S_1^+$ such that $\sqrt{c}u \in \mathcal{G}$ and then
$$
\inf_{u \in \mathcal{G}}E(u) \leq E(\sqrt{c}u)< \inf_{u \in \partial\mathcal{G}}E(u).
$$
Together with Lemma~\ref{EinGandpartialG} we complete the proof. 
\end{proof}

\begin{proof}[\textbf{Proof of Theorem~\ref{theorem1}}]
Let $\nu_c=\inf_{u \in \mathcal{G}}E(u)$ and $\{u_n\}$ be a minimizing sequence of $\nu_c$, i.e. $E(u_n) \to \nu_c$ as $n \to \infty$.
By Lemma~\ref{EinGandpartialG} and Lemma~\ref{EinGandpartialG2}, $\{u_n\}$ is bounded in $H_0^1(\Omega)$, and
 we can assume that $\{u_n\}$ is away from $\partial\mathcal{G}$ passing to a subsequence if necessary. 

By Ekeland's variational principle, we can assume that
$$
(E|_{S_c^+})'(u_n)=(E|_\mathcal{G})'(u_n) \to 0 ~~\text{as}~~ n \to \infty.
$$

Hence, there exists $u_c \in S_c^+$ such that, up to a subsequence, 
\begin{align}\label{untoucweakandstrongandalomost}
  u_n\rightharpoonup  u_c~~\text{ in }H_0^1(\Omega), u_n\to u_c~~\text{ in }L^r(\Omega), \forall~~2\leq r<2^*,
  u_n\to u_c~~\text{a.e. in }\Omega.
\end{align}
We can verify that $u_c \in S_c^+$ is a critical point of $E$ constrained on $S_c^+$.

Let 
$$
\lambda_n=\frac{1}{c}(\int_\Omega f_+(u_n)u_n^+dx-\int_\Omega|\nabla u_n|^{2}dx),
$$
then $\lambda_n$ is bounded and 
\begin{align}\label{Eunandlambdan}
E'(u_n)-\lambda_nu_n^+ \to 0~~\text{in}~~H^{-1}(\Omega)~~\text{as}~~ n \to \infty.
\end{align}

Moreover, there exists $\lambda_c \in \mathbb R$ such that
\begin{align}\label{positivelambdantolambda}
 \lambda_n \to \lambda_c   
\end{align}

and 
\begin{align}\label{Eucandlambdac}
E'(u_c)+\lambda_c u_c=0 ~~\text{in}~~H^{-1}(\Omega).
\end{align}

From~\eqref{untoucweakandstrongandalomost}, \eqref{Eunandlambdan}, \eqref{positivelambdantolambda} and \eqref{Eucandlambdac}, we have
$$
u_n \to  u_c ~~\text{ in }~~H_0^1(\Omega).
$$

As a consequence, we have proved that  $u_c \in \mathcal{G}$  is a critical point of $E|_{S_c^+}$ at the level $\nu_c$.

By lagrange multiplier principle, $u_c$ satisfies
$$
(-\Delta)^s u_c +\lambda_c u_c^+=f_+(u_c)u_c^+~~\text{in}~~ \Omega.
$$

Multiplying $u_c^-$ and integrating on $\Omega$, we obtain $\int_\Omega|\nabla u_c^-|^{2}dx=0$,
which implies that $u_c^-=0$ and hence $u_c \geq 0$.
By the strong maximum principle, $u_c>0$.
Therefore, $\int_\Omega|u_c|^{2}dx=\int_\Omega|u_c^+|^{2}dx=c$ and $(\lambda_c, u_c)$ is a normalized solution of \eqref{Peps}.
Furthermore, 
Assume that $v$ is a normalized solution to ~\eqref{Peps}, then $\int_{\Omega}\vert |v| \vert^2 =\int_{\Omega}\vert v \vert^2=c$. 
Therefore,  
$$
\widetilde{E}(v)=\widetilde{E}(|v|)=E(v)\geq E(u_c)=\widetilde{E}(u_c), 
$$
implying that $u_c$ is a normalized ground state to \eqref{Peps}. 
The proof is completed.
\end{proof}

\section{A mountain pass solution}
In this section, we search for the second positive solution by establishing a mountain pass
structure. 

Define
$$
E_\tau(u):= \frac12\int_\Omega|\nabla u|^2dx-\tau\int_{\Omega}F_+(u)dx
$$
where $\tau \in [\frac{1}{2},1]$.

We can verify that the critical point $u$ of $E_\tau$ on $S_c$ satisfies
\begin{align}\label{equationoftau}
  \begin{cases}-\Delta u+\lambda u=\tau f(u)&\text{in }\Omega,\\
    u>0\quad\text{in }\Omega,\quad u=0&\text{on }\partial\Omega.\end{cases}
\end{align}

\begin{lemma}\label{uandutu}
Let $N \geq 3$, $\Omega$ be bounded, smooth and star-shaped domain with respect to the origin and $f$ satisfy $(f_1)$. 
 If $u\in \mathcal{G}$, then there exists $t=t(u)$ such that $u^{t} \not \in \mathcal{\bar G}$ 
  and $E_\tau(u^t)<0$ uniformly with respect to $\tau$. 
\end{lemma}

\begin{proof}
Recall that $u^t(x)=t^{\frac{N}{2}}u(tx)$. For any $u\in S_c^+$, let

$$
\begin{aligned}
\psi(t)&=\int_\Omega|\nabla u^t|^2dx-(\frac{N}{2}\int_{\Omega}f_+(u)u^+dx-N\int_{\Omega}F_+(u)dx)\\
&\leq \int_\Omega|\nabla u^t|^2dx-\frac{( \alpha-2)N}{2}\int_{\Omega}F_+(u)dx\\
&\leq \int_\Omega|\nabla u^t|^2dx-\frac{( \alpha-2)N\mu}{2} (\int_\Omega| (u^t)^+|^{\alpha}dx+\int_\Omega |(u^t)^+|^{\beta}dx)\\
&=\frac{t^2}{2}\int_\Omega|\nabla u|^2dx-\frac{(\alpha-2)N\mu}{2}(t^{N(\frac{\alpha-2}{2})}\int_\Omega| u^+|^{\alpha}dx
+t^{N(\frac{\beta-2}{2})}\int_\Omega| u^+|^{\beta}dx). 
\end{aligned}
$$
\end{proof}

Moreover, for $\tau \in [\frac{1}{2},1]$, let
$$
  \begin{aligned}
\phi_\tau(t) &=E_\tau(u^t)=\frac12\int_\Omega|\nabla u^t|^2dx-\tau \int_{\Omega}F_+(u^t)dx \\
&\leq \frac12\int_\Omega|\nabla u^t|^2dx-\frac{1}{2} \int_{\Omega}F_+(u^t)dx \\
&\leq \frac12\int_\Omega|\nabla u^t|^2dx-\frac{\mu}{2}  (\int_\Omega| (u^t)^+|^{\alpha}dx+\int_\Omega |(u^t)^+|^{\beta}dx)\\
&=\frac{t^2}{2}\int_\Omega|\nabla u|^2dx-\frac{\mu}{2} \left (t^{N(\frac{\alpha-2}{2})}\int_\Omega| u^+|^{\alpha}dx
+t^{N(\frac{\beta-2}{2})}\int_\Omega| u^+|^{\beta}dx \right ).
    \end{aligned}
$$

Since $N(\frac{\alpha-2}{2}), N(\frac{\beta-2}{2})>2$, we know that
$\psi(t) \to -\infty$ as $t\to +\infty$ and
$\phi_\tau(t) \to -\infty$  uniformly with $\tau$ as $t\to +\infty.$ 
  Therefore, there exists $t$ sufficiently large such that $\psi(t)<0$, implying
  $u^{t} \not \in \mathcal{\bar G}$, and $E_\tau(u^t)=\phi(t)<0$.

\begin{lemma}\label{uniformMP}

  Let $N \geq 3$, $\Omega$ be bounded, smooth and star-shaped domain with respect to the origin, $f$ satisfy $(f_1)$ and
  c satisfy~\eqref{upboundforc}. Then we have
    $$\lim_{\tau \to 1^-} \inf_{u \in \partial \mathcal{G}}E_\tau(u)=\inf_{u \in \partial \mathcal{G}}E(u), $$
 and there exist $\epsilon \in (0,\frac{1}{2})$ and $\delta>0$ such that 
  \begin{align}\label{ucleqpartial}
    E_\tau(u_c)+\delta<\inf_{u \in \partial \mathcal{G}}E_\tau(u), \forall \tau \in (1-\epsilon,1].
  \end{align}
  Moreover, there exists $v\in S_c^+ \backslash \mathcal{G}$ such that
  \begin{align}\label{uniformmtu}
  m_{\tau}:=\inf_{\gamma\in\Gamma}\sup_{t\in[0,1]}E_{\tau}(\gamma(t))>E_{\tau}(u_{c})+\delta=
  \max\{E_{\tau}(u_{c}),E_{\tau}(v)\}+\delta,
  \end{align}
  where 
  $$
  \Gamma=\{\gamma\in C([0,1],S_{c}^{+}):\gamma(0)=u_{c},\gamma(1)=v\}
  $$
  is independent of $\tau$.

\end{lemma}

\begin{proof}

Clearly, $E_\tau(u)\geq E(u)$ for all $u\in \partial \mathcal{G}$ and $\tau \in [\frac{1}{2},1]$, 
implying that
$$
\inf_{u \in \partial \mathcal{G}}E_\tau(u) \geq \inf_{u \in \partial \mathcal{G}}E(u), \forall \tau \in [\frac{1}{2},1].
$$

Let $\{u_n\}$ be a minimizing sequence such that $\lim_{n \to \infty}E(u_n)=\inf_{u\in \partial \mathcal{G}}E(u)$. By \eqref{eugradientu}, 
we know that $\{u_n\}$ is bounded in $H_0^1(\Omega)$. Therefore,
$$
\inf_{u\in \partial \mathcal{G}}E_\tau(u) \leq \liminf_{n \to \infty} E_\tau(u_n)=\lim_{n \to \infty} E(u_n)+o_\tau(1)
=\inf_{u \in \partial \mathcal{G}}E(u)+o_\tau(1),
$$
where $o_\tau(1) \to 0$ as $\tau \to 0^-$. Hence we have

$$
\lim_{\tau \to 1^-}E_\tau(u_c)=E(u_c)=\nu_c. 
$$

By Lemma~\ref{EinGandpartialG2} we have
$$
\lim_{\tau \to 1^-}E_\tau(u_c)=E(u_c)=\nu_c<\inf_{u \in\partial \mathcal{G}}E(u)=\inf_{u \in \partial \mathcal{G}}E_\tau(u).
$$
Choosing $2 \delta=\inf_{u \in \partial \mathcal{G}}E(u)-\nu_c$ and $\epsilon$ sufficiently small we get~\eqref{ucleqpartial}.

By Lemma~\ref{uandutu}, we can let $v=u_c^t$ with $t$ sufficiently large such that $v \not \in \mathcal{\bar G}$ and
$E_\tau(v)<0<E_\tau(u_c)$ for all $\tau \in [\frac{1}{2},1]$.
Since $\gamma(t)$ is continuous, then for any $\gamma \in \Gamma$, 
there exists $t^*\in (0,1]$ such that $\gamma(t^*) \in \partial \mathcal{G}$. Hence 
$$
\inf_{\gamma\in\Gamma}\sup_{t\in[0,1]}E_{\tau}(\gamma(t))\geq \inf_{u \in \partial \mathcal{G}}E_\tau(u).
$$
Therefore, we obtain~\eqref{uniformmtu}. 

\end{proof}

\begin{lemma}\label{matautom1}
  Under the assumptions of Lemma~\ref{uniformMP}, we have $\lim_{\tau \to 1^-}m_\tau=m_1$.
\end{lemma}

It follows from the definition of $E_\tau$ and $m_\tau$ that $m_\tau \geq m_1$ for any $\tau \in [\frac{1}{2},1)$, 
then $\liminf_{\tau \to 1^-}m_\tau \geq m_1$.
It is sufficient to prove $\limsup_{\tau \to 1^-}m_\tau \leq m_1$.
For any $\epsilon>0$, we can take $\gamma_0$ such that
$$
\sup_{t\in [0,1]}E(\gamma_0(t))<m_1+\epsilon.
$$
For any $\tau_n \to 1^-$, we have
$$
m_{\tau_n}=\inf_{\gamma \in \Gamma}\sup_{t\in [0,1]}E_{\tau_n}(\gamma(t)) \leq \sup_{t\in [0,1]}E_{\tau_n}(\gamma_0(t))
=\sup_{t\in [0,1]}E(\gamma_0(t))+o_n(1)<m_1+\epsilon+o_n(1).
$$
By the arbitrariness of $\epsilon$ we obtain that
$$
\lim_{n \to \infty}m_{\tau_n} \leq m_1,
$$
implying that
$$
\limsup_{\tau \to 1^-}m_{\tau_n} \leq m_1.
$$
Therefore, 
$$
\lim_{\tau \to 1^-}=m_1.
$$

We introduce the monotonicity trick from~\cite{jeanjean1999existence} on the family of functionals $E_\tau$ to obtain a bounded (PS) sequence. 
Indeed, we use the version applied to $H_0^1(\Omega)$ and $S_c^+$ as follows. 

\begin{theorem}\label{Monotonicitytrick} (Monotonicity trick).

Let $I=[\frac{1}{2},1]$.  
We consider a family $(E_\tau)_{\tau\in I}$ of $C^1$-functionals on $H_0^1(\Omega)$ of the form
$$
E_\tau(u)=A(u)-\tau B(u),\quad\tau\in I, 
$$
where $B(u)\geq0,\forall u\in H_0^1(\Omega)$ and such that either $A( u) \to \infty~~\textit{or }~~B( u) \to \infty$ as $\| u\| \to \infty.$ 
We assume there are two points $(u_1,u_2)$ in $S_c^+$ (independent of $\tau$) such that setting
 $$
\Gamma=\{\gamma\in C([0,1],S_c^+),\gamma(0)=u_1,\gamma(1)=u_2\},
$$
there  holds, $\forall \tau \in I$,
$$
m_\tau:=\inf_{\gamma\in\Gamma}\sup_{t\in[0,1]}E_\tau(\gamma(t))>\max\{E_\tau(u_1),E_\tau(u_2)\}.
$$
Then, for almost everywhere $\tau\in I$, there is a sequence $\{u_n\}\subset S_c^+$ such that
\begin{itemize}
  \item[$(i)$] $\{u_n\}$ is bounded in $H_0^1(\Omega)$;
  \item[$(ii)$] $E_\tau ( u_n) \to m_\tau$;
  \item[$(iii)$] $E_\tau ^{\prime }| _{S_c^+}( u_n) \to 0$ in $H^{-1}(\Omega)$.
  \end{itemize} 

\end{theorem}

\begin{proof}
Referring to~\cite[Lemma 3.5]{jeanjean1999existence} and~\cite{Qishijie}, we provide a concise proof framework. 
Denote $m_\tau$ by $m(\tau)$. 
It follows from the definition of $E_\tau$ that $m(\tau) \geq m(1)$ for any $\tau \in [\frac{1}{2},1]$ and 
$m(\tau)$ is nonincreasing with respect to $\tau$. 
By Lebesgue Theorem, $m(\tau)$ is a.e. differentiable at $[\frac{1}{2},1]$. 
Fix $\tau\in [\frac{1}{2},1]$ such that $m'(\tau)$ exist from now on. 
Let $\{\tau_n\}$ be a sequence such that $\tau_n <\tau$ and $\tau \to \tau$ as $n \to \infty$. 
Then there exist $n(\tau)>0$ such that 
\begin{align}\label{Etaugamman0}
  -m'(\tau)-1<\frac{m(\tau)-m(\tau_n)}{\tau-\tau_n}<-m'(\tau)+1,~~\forall n \geq n(\tau). 
\end{align}

First, we prove that if we impose certain limitations on the oscillation range of $\gamma_n$, 
there is an upper bound for $|A(\gamma_n(t))|$ only with respect to $m'(\tau)$.
Specifically, we prove that there exist $\gamma_n \in \Gamma$ and $M=M(m'(\tau))$ such that 
if 
\begin{align}\label{Etaugamman01}
E_\tau(\gamma_n(t))\geq m(\tau)-(\tau-\tau_n)
\end{align}

for some $t \in [0,1]$, then $|A(\gamma_n(t))| \leq M$. 
Indeed, letting $\gamma_n \in \Gamma$ be such that
\begin{align}\label{Etaugamman1}
\sup_{t\in[0,1]}E_{\tau_{n}}(\gamma_{n}(t))\leq m(\tau_{n})+(\tau-\tau_{n}), 
\end{align}
we have
\begin{align}\label{Etaugamman2}
  E_\tau(\gamma_n(t))\leq E_{\tau_n}(\gamma_n(u))\leq m(\tau_n)+(\tau-\tau_n)\leq m(\tau)+(2-m^{\prime}(\tau))(\tau-\tau_n).  
\end{align}

Moreover, if $E_\tau(\gamma_n(t))\geq m(\tau)-(\tau-\tau_n)$ for some $t\in[0,1]$,

then
\begin{align}\label{Etaugamman3}
E_{\tau_n}(\gamma_n(t))-E_\tau(\gamma_n(t))=(\tau-\tau_n)B(\gamma_n(t)).
\end{align}

On the other hand, by~\eqref{Etaugamman01} and~\eqref{Etaugamman2}, we have
\begin{align}\label{Etaugamman4}
E_{\tau_n}(\gamma_n(t))-E_\tau(\gamma_n(t))\leq m(\tau_n)-m(\tau)+2(\tau-\tau_n). 
\end{align}

Combining~\eqref{Etaugamman2}, ~\eqref{Etaugamman3} and~\eqref{Etaugamman4}, we have
$$
B(\gamma_n(t)) \leq 3-m^{\prime}(\tau). 
$$

$$
  \begin{aligned}
    |A(\gamma_n(t))| & \leq E_{\tau_n}(\gamma_n(t))+\tau_nB(\gamma_n(t)) \\
    &\leq m(\tau_n)+(\tau-\tau_n)+\tau_n(3-m^{\prime}(\tau))\\
    &\leq m(\tau)+(1-m^{\prime}(\tau))(\tau-\tau_n)+(\tau-\tau_n)+\tau_n(3-m^{\prime}(\tau))\\
    &\leq m(\tau)+3\tau-m^{\prime}(\tau)\tau \\
    &\leq m(\frac{1}{2})+\tau(3-m^{\prime}(\tau)).
    \end{aligned}
$$

Secondly, we prove that when $E_\tau(u)$ approaches to $m(\tau)$ with
 $|A(u)|$ controlled by a constant depending on $M$, there exist a bounded Palais-Smale sequence at the level $m(\tau)$. That is,  
similar to the proof in~\cite[Lemma 3.5]{Qishijie}, we can establish a deformation on $S_c^+$. 
Moreover, setting 
$$
B_1(\epsilon)=\{u \in S_c^+: |E_\tau-m(\tau)| < \epsilon, |A(u)| \leq M+1\},
$$
where $\epsilon=(m(\tau)-\max\{E_\tau(u_1),E_\tau(u_2)\})$.
 we obtain a bounded Palais-Smale sequence at the level $m(\tau)$ in $B_1(\epsilon)$.

\end{proof}

\begin{proposition}\label{criticalpointatetau}
Under the assumptions of Lemma~\ref{uniformMP}, 
for almost everywhere $\tau \in [1-\epsilon, 1]$ where $\epsilon$ is given by Lemma~\ref{uniformMP}, there exists a critical point $u_\tau$
of $E_\tau$ constrained on $S_c^+$ at the level $m_\tau$, which solves \eqref{equationoftau}.

\end{proposition}

\begin{proof}

Let
$$
A(u)=\frac{1}{2}\int_\Omega|\nabla u|^{2}dx,~~ B(u)=\int_\Omega F_+(u)dx.
$$

By Theorem~\ref{Monotonicitytrick} and Lemma~\ref{uniformMP}, 
we know that for almost everywhere $\tau \in [1-\epsilon,1]$, there exists a bounded (PS) sequence $\{u_n\} \subset S_c^+$ 
satisfying $E_\tau(u_n) \to m_\tau$ and $(E_\tau|_{S_c^+})'(u_n) \to 0$ as $n \to \infty$.
Therefore, there exist $u_c \in S_c^+$ such that, up to a subsequence, 
\begin{align}\label{untoutauweakandstrongandalomost}
  u_n\rightharpoonup  u_\tau~~\text{ in }H_0^1(\Omega), u_n\to u_\tau~~\text{ in }L^r(\Omega), \forall2\leq r<2^*,
  u_n\to u_\tau~~\text{a.e. in }\Omega.
\end{align}

We can verify that $u_\tau \in S_c^+$ is a critical point of $E_\tau$ constrained on $S_c^+$.

Let 
$$
\lambda_n=\frac{1}{c}(\tau \int_\Omega  f(u_n^+)u_n^+dx-\int_\Omega|\nabla u_n|^{2}dx),
$$
then $\lambda_n$ is bounded and 
\begin{align}\label{Etauunandlambda}
E_\tau'(u_n)-\lambda_nu_n^+ \to 0~~\text{in}~~H^{-1}(\Omega)~~\text{as}~~ n \to \infty.
\end{align}

Moreover,  there exist $\lambda_\tau \in \mathbb R$ such that
\begin{align}\label{lambdantotau}
\lambda_n \to \lambda_\tau
\end{align}
and 
\begin{align}\label{Etauutauandlambda}
  E'(u_\tau)+\lambda_\tau u_\tau=0 ~~\text{in}~~H^{-1}(\Omega).
\end{align}.

From~\eqref{untoutauweakandstrongandalomost}, \eqref{lambdantotau},  \eqref{Etauunandlambda} and~\eqref{Etauutauandlambda}, we have
$$
u_n \to  u_\tau ~~\text{ in }~~H_0^1(\Omega).
$$

By lagrange multiplier principle, $u_\tau$ satisfies
$$
(-\Delta)^s u_\tau +\lambda_\tau u_\tau^+=\tau f_+(u_\tau)u_\tau^+~~\text{in}~~ \Omega.
$$

Multiplying $u_\tau^-$ and integrating on $\Omega$, we obtain $\int_\Omega|\nabla u_\tau^-|^{2}dx=0$,
which implies that $u_\tau^-=0$ and hence $u_\tau \geq 0$.
By the strong maximum principle, $u_\tau>0$.
Therefore, $\int_\Omega|u_\tau|^{2}dx=\int_\Omega|u_\tau^+|^{2}dx=c$ and $(\lambda_\tau, u_\tau)$ is a normalized solution of \eqref{equationoftau}.
\end{proof}

\begin{proof}[\textbf{Proof of Theorem~\ref{theorem2}}]
  By Lemma~\ref{matautom1}, we can take $\epsilon$ sufficiently small such that 
$$
 \forall \tau \in [1-\epsilon,1], m_{\tau} \leq 2m_1.
$$

Let $\{\tau_n\} \subset [1-\epsilon,1]$ be a sequence  such that $\tau_n \to 1^-$ as $n \to \infty$, and $(\lambda_n, u_n)$ is a
normalized solution of~\eqref{equationoftau} at the level $m_{\tau_n}$. Then we have
$$
2m_1\geq m_{\tau_n}=E_{\tau_n}(u_n)=\frac{1}{2}\int_\Omega|\nabla u_n|^2dx-\tau_n \int_{\Omega}F_+(u_n)dx
\geq (\frac{1}{2}-\frac{2\tau_n}{(\alpha-2) N})\int_\Omega|\nabla u_n|^2dx, 
$$
hence $\{u_n\} $ is bounded in $H_0^1(\Omega)$. 

It follows from ~\eqref{equationoftau} that
$$
\lambda_n=\frac{1}{c}(\tau_n\int_\Omega f(u_n)u_ndx-\int_\Omega|\nabla u|^{2}dx),
$$
and so $\lambda_n$ is bounded.
Therefore, there exists $\widetilde{u}_c$ and $\widetilde{\lambda}_c$ such that
\begin{align}\label{widetildeuntoutauweakandstrongandalomost}
  u_n\rightharpoonup  \widetilde{u}_c~~\text{ in }H_0^1(\Omega), u_n\to \widetilde{u}_c~~\text{ in }L^r(\Omega), \forall2\leq r<2^*,
  u_n\to \widetilde{u}_c~~\text{a.e. in }\Omega,
\end{align}
and
\begin{align}\label{lambdantowidec}
\lambda_n \to \widetilde{\lambda}_c. 
\end{align}

By~\eqref{widetildeuntoutauweakandstrongandalomost}, \eqref{lambdantowidec} and~~\eqref{equationoftau}, we have
$$
u_n \to \widetilde{u}_c~~\text{in}~~H_0^1(\Omega).
$$
Consequently,  $\int_\Omega|\widetilde{u}_c|^{2}dx=c$ and $\widetilde{u}_c \geq 0$.
By the strong maximum principle, $u_c>0$.
Therefore, $E(\widetilde{u}_c)=m_1$, and $\widetilde{u}_c$ is 
a positive normalized solution of \eqref{Peps}.

\end{proof}

\bigskip

\section{some mixed cases}
In this section, we consider~\eqref{Peps} when $f$ is mixed with $L^2$-supercritical terms and critical or subcritical terms.
We consider the case that
  \begin{align}\label{formulap+q}
    f(u)= |u|^{p-2}u+a|u|^{q-2}u,~2<q<2+\frac{4}{N}<p<2^*, 
     \end{align}
then
\begin{align}
  E(u)=\frac12\int_\Omega|\nabla u|^2dx-\frac{1}{p}\int_\Omega|u^+|^p dx-\frac{a}{q}\int_\Omega|u^+|^q dx,
\end{align}

$$
\mathcal{G}:=\left\{u\in S_c^+:\int_\Omega|\nabla u|^2dx>
(\frac{1}{2}-\frac{1}{p})N\int_\Omega|u^+|^p dx+a(\frac{1}{2}-\frac{1}{q})N\int_\Omega|u^+|^qdx\right\}.
$$
and
$$
\partial \mathcal{G}:=\left\{u\in S_c^+:\int_\Omega|\nabla u|^2dx=
(\frac{1}{2}-\frac{1}{p})N\int_\Omega|u^+|^p dx +a(\frac{1}{2}-\frac{1}{q})N\int_\Omega|u^+|^q dx\right\}.
$$

For the convenience of the following arguments, first, we introduce the Gagliardo-Nirenberg inequality. 
 For any $N\geq2$ and $p\in(2,2^*)$, there is a constant $C_{N,p}$ depending on
  $N$ and $p$ such that
  $$
  \int_{\mathbb{R}^N}|u|^pdx
  \leq C_{N,p}^p\left(\int_{\mathbb{R}^N}|u|^2dx\right)^{\frac{2p-N(p-2)}4}\left(\int_{\mathbb{R}^N}|\nabla u|^2dx\right)^{\frac{N(p-2)}4},\:\forall\:u\in H^1(\mathbb{R}^N),
  $$
  where the optimal constant $C_{N,p}$ can be expressed exactly as
  $$
  C_{N,p}^p=\frac{2p}{2N+(2-N)p}\left(\frac{2N+(2-N)p}{N(p-2)}\right)^{\frac{N(p-2)}4}\frac1{\left\|Q_p\right\|_2^{p-2}},
  $$
  and $Q_p$ is the unique positive solution, up to translations, of the equation
  $$-\Delta Q+Q=|Q|^{p-2}Q\quad\text{in }\mathbb{R}^N.$$

It is sufficient to prove the following lemma to obtain our conclusions. 

\begin{lemma}\label{anotherEinGandpartialG}
 Let $N \geq 3$, $\Omega$ be bounded, smooth and star-shaped domain with respect to the origin, $f=|u|^{p-2}u+a|u|^{q-2}u$, and $c$ satisfy
 \eqref{upuqccage0} and \eqref{upuqccale0}.
 Then 
$\inf_{u \in \mathcal{G}}E(u)$ has a lower bound, and
any sequence $\{u_n\} \subset \mathcal{G}$ satisfying $\limsup_{n\to \infty} E(u_n) < +\infty$ is bounded in $H_0^1(\Omega)$. 
Moreover, $\mathcal{G} \neq \varnothing$, and 
there holds that
$$
\inf_{u\in \mathcal{G}} E(u)<\inf_{u\in \partial \mathcal{G}}E(u).
$$

   \end{lemma}

\begin{proof}
  
For any $u \in \mathcal{G}$, direct calculation yields that
\begin{align}\label{aformulaofupuq}
  \frac{2}{N(p-2)}\int_\Omega|\nabla u|^2dx-\frac{a(q-2)}{q(p-2)}\int_\Omega|u|^q dx>\frac{1}{p}\int_\Omega|u|^p dx.
\end{align}

First, we consider the case that $a>0$. 
By~\eqref{aformulaofupuq}, we have
$$
  \begin{aligned}
    E(u) &>(\frac{1}{2}-\frac{2}{N(p-2)})\int_\Omega|\nabla u|^2dx-\frac{a(p-q)}{q(p-2)}\int_\Omega|u^+|^q dx\\
    &\geq (\frac{1}{2}-\frac{2}{N(p-2)})\int_\Omega|\nabla u|^2dx-\frac{a(p-q)}{q(p-2)}C_{N,q}^q c^{\frac{2q-N(q-2)}4}
    \left(\int_{\mathbb{R}^N}|\nabla u|^2dx\right)^{\frac{N(q-2)}4}.
    \end{aligned}
$$
 
Since $0<\frac{N(q-2)}4<1$, 
recalling that 
\begin{align}\label{nablaugeqlambda1cagain}
  \int_\Omega|\nabla u|^2dx \geq \lambda_1c>0, 
 \end{align}
we can verify that $E(u)$ is bounded from below on $\mathcal{G}$. 

Let $\{u_n\} \subset \mathcal{G}$ satisfy $\limsup_{n\to \infty} E(u_n) < +\infty$. 
Assume by contradiction that there exist a subsequence of $\{u_n\}$ (still denoted by $\{u_n\}$) 
such that $\lim_{n \to \infty}\int_\Omega|\nabla u|^2dx \to +\infty$ as $n \to \infty$, 
then 
$$
(\frac{1}{2}-\frac{2}{N(p-2)})\int_\Omega|\nabla u|^2dx-\frac{a(p-q)}{q(p-2)}C_{N,q}^q c^{\frac{2q-N(q-2)}4}
\left(\int_{\mathbb{R}^N}|\nabla u|^2dx\right)^{\frac{N(q-2)}4} \to +\infty,
$$
as $n \to \infty$, implying that $E(u_n) \to +\infty$ as $n \to \infty$. Contradiction! 

For any $u\in S_1^+$ and $c$ satisfying $(\frac{1}{2}-\frac{1}{p})c^{\frac{p-2}{2}}+a(\frac{1}{2}-\frac{1}{q})c^{\frac{q-2}{2}}<g_{3,u}$, 
where $$g_{3,u}=\frac{\int_{\Omega}|\nabla u|^2dx}{\int_{\Omega}|u^+|^{p}dx+a\int_{\Omega}|u^+|^{q}dx}, $$
there holds that
$$
  \begin{aligned}
    \int_\Omega|\nabla (\sqrt{c}u)|^{2}dx &=c\int_\Omega|\nabla u|^{2}dx \\
    &>(\frac{1}{2}-\frac{1}{p})c^{\frac{p}{2}}\int_\Omega|u^+|^{p}dx+a(\frac{1}{2}-\frac{1}{q})c^{\frac{q}{2}}\int_\Omega|u^+|^{q}dx\\
    &=(\frac{1}{2}-\frac{1}{p})\int_\Omega| (\sqrt{c}u^+)|^{p}dx+a(\frac{1}{2}-\frac{1}{q})\int_\Omega| (\sqrt{c}u^+)|^{q}dx.
    \end{aligned}
$$ 

Hence $\sqrt{c}u \in \mathcal{G}$ and $\mathcal{G}$ is not empty.

For any $u\in \partial \mathcal{G}$, by~\eqref{sobinequality}, we have
$$
  \begin{aligned}
    \int_\Omega|\nabla u|^2dx &=(\frac{1}{2}-\frac{1}{p})N\int_\Omega|u^+|^p dx+a(\frac{1}{2}-\frac{1}{q})N\int_\Omega|u^+|^q dx\\
    &\leq (\frac{1}{2}-\frac{1}{p})(S_{p^*}^{-1}|\Omega|^{\frac{2-p^*}{2}})^{\frac{p}{p^*}} (\int_\Omega|\nabla u|^2dx)^{\frac{p}{2}}\\
    & \qquad +a(\frac{1}{2}-\frac{1}{q})N(S_{q^*}^{-1}|\Omega|^{\frac{2-q^*}{2}})^{\frac{q}{q^*}} (\int_\Omega|\nabla u|^2dx)^{\frac{q}{2}}. 
    \end{aligned}
$$   
Since $\frac{q}{2}, \frac{p}{2}>1$, together with \eqref{nablaugeqlambda1cagain}, there is $C_2=C_2(p,q,N)>0$ such that  
\begin{align}\label{nablageqC2}
  \int_\Omega|\nabla u|^2dx \geq C_2. 
\end{align}
Indeed, $C_2$ satisfies
\begin{align}\label{defineC2}
(\frac{1}{2}-\frac{1}{p})(S_{p^*}^{-1}|\Omega|^{\frac{2-p^*}{2}})^{\frac{p}{p^*}} (C_2)^{\frac{p}{2}}
  +a(\frac{1}{2}-\frac{1}{q})N(S_{q^*}^{-1}|\Omega|^{\frac{2-q^*}{2}})^{\frac{q}{q^*}} (C_2)^{\frac{q}{2}}=1. 
\end{align}

It is not difficult to verify that $C_2>0$ satisfying~\eqref{defineC2} is unique.

Therefore, for any $c$ satisfying
$$
\frac{a(p-q)}{q(p-2)}C_{N,q}^q c^{\frac{2q-N(q-2)}4}
<(\frac{1}{4}-\frac{1}{N(p-2)})C_2^{\frac{4-N(q-2)}4}, 
$$
and $u \in \partial\mathcal{G} $, 
we have
$$
  \begin{aligned}
    E(u) &\geq (\frac{1}{2}-\frac{2}{N(p-2)})\int_\Omega|\nabla u|^2dx-\frac{a(p-q)}{q(p-2)}C_{N,q}^q c^{\frac{2q-N(q-2)}4}
    (\int_{\mathbb{R}^N}|\nabla u|^2dx)^{\frac{N(q-2)}4} \\
    &=\left( (\frac{1}{2}-\frac{2}{N(p-2)})-\frac{a(p-q)}{q(p-2)}C_{N,q}^q c^{\frac{2q-N(q-2)}4}
    (\int_{\mathbb{R}^N}|\nabla u|^2dx)^{\frac{N(q-2)-4}4} \right)\int_\Omega|\nabla u|^2dx\\
    &> (\frac{1}{4}-\frac{1}{N(p-2)})C_2.
    \end{aligned}
$$

For any $u \in S_c^+$ satisfying $c<(\frac{1}{4}-\frac{1}{N(p-2)})C_2g_{2,u}^{-1}$, 
we have
$$
  \begin{aligned}
    E(\sqrt{c}u) &
    =\frac{c}{2}\int_\Omega|\nabla u|^{2}dx-\frac{c^\frac{p}{2}}{p}\int_\Omega|u^+|^{p}dx-\frac{ac^\frac{q}{2}}{q}\int_\Omega|u^+|^{q}dx \\
    &<\frac{c}{2}\int_\Omega|\nabla u|^{2}dx\\
    &<(\frac{1}{4}-\frac{1}{N(p-2)})C_2.
    \end{aligned}
$$

When ~\eqref{upuqccage0} holds, we can take $u\in S_1^+$ such that $\sqrt{c}u \in \mathcal{G}$ and then
$$
\inf_{u \in \mathcal{G}}E(u) \leq E(\sqrt{c}u)< \inf_{u \in \partial\mathcal{G}}E(u).
$$

Secondly, we consider the case that $a<0$. 
By~\eqref{aformulaofupuq} and \eqref{nablaugeqlambda1cagain}, we can deduce that 
$$
  \begin{aligned}
    E(u) &>(\frac{1}{2}-\frac{2}{N(p-2)})\int_\Omega|\nabla u|^2dx-\frac{a(p-q)}{q(p-2)}\int_\Omega|u^+|^q dx\\
    &\geq (\frac{1}{2}-\frac{2}{N(p-2)})\int_\Omega|\nabla u|^2dx\\
    &>0, 
    \end{aligned}
$$
then
$$
\limsup_{n \to \infty}(\frac12-\frac{2}{N( p -2)})\int_\Omega|\nabla u|^2dx<\limsup_{n \to \infty}E(u_n)<+\infty. 
$$
Therefore, $\{u_n\} \subset \mathcal{G}$ is bounded in $H_0^1(\Omega)$.

For any $u\in S_1^+$, $c$ satisfying $(\frac{1}{2}-\frac{1}{p})c^{\frac{p-2}{2}}<g_{4,u}$, 
where
$$
g_{4,u}=\frac{\int_{\Omega}|\nabla u|^2dx}{\int_{\Omega}|u^+|^{p}dx}, 
$$
we have
$$
  \begin{aligned}
    \int_\Omega|\nabla (\sqrt{c}u)|^{2}dx &=c\int_\Omega|\nabla u|^{2}dx \\
    &>(\frac{1}{2}-\frac{1}{p})c^{\frac{p}{2}}\int_\Omega|u^+|^{p}dx\\
    &>(\frac{1}{2}-\frac{1}{p})\int_\Omega| (\sqrt{c}u^+)|^{p}dx+a(\frac{1}{2}-\frac{1}{q})\int_\Omega| (\sqrt{c}u^+)|^{q}dx.
    \end{aligned}
$$ 

Hence $\sqrt{c}u \in \mathcal{G}$ and $\mathcal{G}$ is not empty.

For any $u\in \partial \mathcal{G}$, by~\eqref{sobinequality}, we have
$$
  \begin{aligned}
    \int_\Omega|\nabla u|^2dx &=(\frac{1}{2}-\frac{1}{p})N\int_\Omega|u^+|^p dx+a(\frac{1}{2}-\frac{1}{q})N\int_\Omega|u^+|^q dx\\
    &\leq (\frac{1}{2}-\frac{1}{p})(S_{p^*}^{-1}|\Omega|^{\frac{2-p^*}{2}})^{\frac{p}{p^*}} (\int_\Omega|\nabla u|^2dx)^{\frac{p}{2}}
    \end{aligned}
$$

Since $\frac{q}{2}, \frac{p}{2}>1$, there is $C_3=C_3(p,q,N)>0$ such that  
\begin{align}\label{nablageqC3}
  \int_\Omega|\nabla u|^2dx \geq C_3, 
\end{align}
Indeed, $C_3$ satisfies
\begin{align}\label{defineC3}
  (\frac{1}{2}-\frac{1}{p})(S_{p^*}^{-1}|\Omega|^{\frac{2-p^*}{2}})^{\frac{p}{p^*}} (C_3)^{\frac{p-2}{2}}=1. 
\end{align}
It is not difficult to verify that $C_3>0$ satisfying~\eqref{defineC3} is unique.

Therefore, 
there holds that
$$
\inf_{u \in \partial\mathcal{G}}E(u)\geq (\frac{1}{2}-\frac{2}{N(p-2)})\int_\Omega|\nabla u|^2dx
\geq (\frac{1}{2}-\frac{2}{N(p-2)})C_3.
$$

For any $u\in S_1^+$ and 
$$
c<\max_{\varsigma \in (0,1)} \left\{min\left\{\varsigma(\frac{1}{2}-\frac{2}{N(p-2)})C_3g_{2,u}^{-1},~
\left ((1-\varsigma)\frac{q}{|a|}(\frac{1}{2}-\frac{2}{N(p-2)})C_3g_{5,u}\right )^{\frac{2}{q}}\right\}\right\},
$$
where
$$
g_{5,u}=\frac{1}{\int_{\Omega}|u^+|^{q}dx}, 
$$

there exists a 
\begin{align}\label{definevarsigma0}
  \varsigma_0 \in (0,1)
\end{align}

such that 
$$
  \begin{aligned}
    E(\sqrt{c}u) &
    =\frac{c}{2}\int_\Omega|\nabla u|^{2}dx-\frac{c^\frac{p}{2}}{p}\int_\Omega|u^+|^{p}dx-a\frac{c^\frac{q}{2}}{q}\int_\Omega|u^+|^{q}dx \\
    &<\frac{c}{2}\int_\Omega|\nabla u|^{2}dx+|a|\frac{c^\frac{q}{2}}{q}\int_\Omega|u^+|^{q}dx\\
    &<\varsigma_0(\frac{1}{2}-\frac{2}{N(p-2)})C_3+(1-\varsigma_0)(\frac{1}{2}-\frac{2}{N(p-2)})C_3\\
    &<(\frac{1}{2}-\frac{2}{N(p-2)})C_3.
    \end{aligned}
$$

we can take $u\in S_1^+$ such that $\sqrt{c}u \in \mathcal{G}$ and then
$$
\inf_{u \in \mathcal{G}}E(u) \leq E(\sqrt{c}u)< \inf_{u \in \partial\mathcal{G}}E(u).
$$ 

Together with~\eqref{nablaugeqlambda1cagain}, we complete the proof. 

\end{proof}

\begin{proof}[\textbf{Proof of Theorem~\ref{theorem3}}]
Similar to the proof of Theorem~\ref{theorem1} and Theorem~\ref{theorem2}, we can obtain our results. 
\end{proof}

\section{Multiple normalized solutions}

The author in~\cite{Qishijie} obtain the existence of multiple normalized solutions of~\eqref{Peps} when $f=|u|^{p-2}u+|u|^{q-2}u$,  
where  $2<q<2+\frac{4}{N}<p<2^*$. In this section, 
we consider the multiplicity of solutions of~\eqref{Peps} 
with either of the following conditions:
\begin{itemize}
  \item [(1)]$f$ satisfies $(f_1)$; 
  \item [(2)]$f=|u|^{p-2}u-|u|^{q-2}u$, where  $2<q<p$, $2+\frac{4}{N}<p<2^*$. 
\end{itemize}

Let $0<\theta_1<\theta_2< \cdots< \theta_i<\cdots$ be the sequence of different Dirichlet eigenvalues of $-\Delta$ on $\Omega$, 
$\delta_i$ be the multiplicity of $\theta_i$, and $e_{i,j}(j=1,2,\cdots, \delta_i)$ be the corresponding orthonormal eigenfunctions 
in $L^2(\Omega)$.
Define $V_i=span\{e_{1}, \cdots, e_{i1}, \cdots, e_{i\delta_i}\}$,
then 
$$
V_1 \subset \cdots \subset V_i \subset V_{i+1} \subset \cdots,~~\text{ and }~~\cup_{i=1}^{+\infty}=H_0^1(\Omega). 
$$

\subsection{Multiple normalized solutions with (1) holds}

For any $\tau \in [\frac{1}{2},1]$, we define $\widetilde{E}_\tau(u):S_c \to \mathbb R^N$ by 
$$
\widetilde{E}_\tau(u) = \frac12\int_\Omega|\nabla u|^2dx-\tau \int_\Omega F(u)dx. 
$$

Clearly, the critical point $u$ of $\widetilde{E}_\tau(u)$ on $S_c$ satisfies
\begin{align}\label{fequationoftauwidetildeE}
  \begin{cases}-\Delta u+\lambda u=\tau f(u)&\text{in }\Omega, \quad u_{|\partial \Omega}=0, \\
    \int_{\Omega}\vert u \vert^2 =c.\end{cases}
\end{align}

Define $\rho_i$ by
\begin{align}\label{definerhoi}
  \zeta (C_{N,\alpha}^\alpha \theta_i^{{\frac{N(\alpha-2)-2\alpha}4}}\rho_i^{\frac{\alpha-2}2}
  +C_{N,\beta}^\beta \theta_i^{{\frac{N(\beta-2)-2\beta}4}}\rho_i^{\frac{\beta-2}2})=\frac{1}{2}.
\end{align}

Since $2<\alpha<\beta<2^*$, then $\frac{\alpha-2}2, \frac{\beta-2}2>0$ and $\frac{N(\alpha-2)-2\alpha}4, \frac{N(\beta-2)-2\beta}4<0$.
Therefore, $\rho_i$ is well defined and $\rho_i \to +\infty$ as $i \to \infty$.

\begin{lemma}\label{lemmaEtauvonBigeqEtauu}
  Let
  \begin{align}\label{defineci}
    \widetilde{c}_i=\frac{\rho_i}{2\theta_i}
  \end{align}
and  $c<\widetilde{c}_i$. Define
$$
\mathcal{B}_i=\{v \in V_{i-1}^\bot \cap S_c:\|\nabla v\|_2^2=\rho_i \}.
$$ 
Then for any $u\in V_i \cap S_c$, we have
$$
\int_\Omega|\nabla u|^{2}dx<\rho_i, 
$$
and 
\begin{align}\label{EtauvonBigeqEtauu}
  \widetilde{E}_\tau(u)< \inf_{v \in \mathcal{B}_i}\widetilde{E}_\tau(v) ~~\text{for any}~~\tau \in[\frac{1}{2},1].
\end{align}

\end{lemma}

\begin{proof}

For any $ u \in V_{i} \cap S_c$, we have $\int_\Omega|\nabla u|^2dx\leq \theta_i c<\frac{\rho_i}{2}$. Then
$$
\widetilde{E}_\tau(u) \leq \frac12\int_\Omega|\nabla u|^2dx - \int_\Omega F(u) dx \leq \frac{1}{2}\theta_i c
< \frac{1}{4}\rho_i. 
$$

On the other hand, for any $v\in \mathcal{B}_i$,

$$
  \begin{aligned}
    \widetilde{E}_\tau(v) &\geq \frac12\int_\Omega|\nabla v|^2dx-\frac{1}{2}\zeta (\int_\Omega|v|^\alpha dx +\int_\Omega|v|^\beta dx)\\
  &\geq \frac{1}{2}\rho_i-\frac{1}{2}\zeta (C_{N,\alpha}^\alpha c^{\frac{2\alpha-N(p-2)}4}\rho_i^{\frac{N(\alpha-2)}4}
  +C_{N,\beta}^\beta c^{\frac{2\beta-N(p-2)}4}\rho_i^{\frac{N(\beta-2)}4} )\\
  &= \left(\frac{1}{2}-\frac{1}{2}\zeta (C_{N,\alpha}^\alpha \theta_i^{{\frac{N(\alpha-2)-2\alpha}4}}\rho_i^{\frac{\alpha-2}2}
  +C_{N,\beta}^\beta \theta_i^{{\frac{N(\beta-2)-2\beta}4}}\rho_i^{\frac{\beta-2}2})\right)\rho_i
  \end{aligned}
$$

It follows from the definition of $\rho_i$ that $\widetilde{E}_\tau(u)< \inf_{v \in \mathcal{B}_i}\widetilde{E}_\tau$.

\end{proof}

\begin{lemma}\label{findwideu}
  For any $u\in V_i \cap S_c$ and $t>0$, define
\begin{align}\label{definewideu}
  u^t(x):=t^{\frac{N}{2}}u(tx). 
\end{align}
Then there exists $t=t(i,c)$ such that $\int_\Omega|u^t|^2=\int_\Omega|u|^2=c$,
$\int_\Omega|\nabla u^t|^2dx \geq 2\rho_i$ and
$\widetilde{E}_\tau (u^t)< 0$ uniformly with respect to $\tau$ and $ u\in V_i \cap S_c$.

\end{lemma}

\begin{proof}

We can verify that $\int_\Omega|u^t|^2=\int_\Omega|u|^2=c$,  
$\int_\Omega|\nabla u^t|^2dx=t^2\int_\Omega|\nabla u|^2dx$,
and 
$$
  \begin{aligned}
\widetilde{E}_\tau(u^t) &=\frac12\int_\Omega|\nabla u^t|^2dx-\tau \int_{\Omega}F_+(u^t)dx \\
&\leq \frac{t^2}{2}\int_\Omega|\nabla u|^2dx-\frac{\mu}{2} \left (t^{N(\frac{\alpha-2}{2})}\int_\Omega| u|^{\alpha}dx
+t^{N(\frac{\beta-2}{2})}\int_\Omega| u|^{\beta}dx \right )\\
& \leq \frac{t^2}{2} \theta_i c-\frac{\mu}{2} \left (t^{N(\frac{\alpha-2}{2})} c^{\frac{\alpha}{2}}|\Omega|^{\frac{2-\alpha}{2}}
+t^{N(\frac{\beta-2}{2})}c^{\frac{\beta}{2}}|\Omega|^{\frac{2-\beta}{2}} \right ).
    \end{aligned}
$$
Since $N(\frac{\alpha-2}{2}), N(\frac{\beta-2}{2})>2$, 
then $\widetilde{E}_\tau(u^t) \to -\infty$ as $t \to +\infty$ uniformly with respect to $\tau \in [\frac{1}{2},1]$ and $ u\in V_i \cap S_c$. 

Therefore, we can take $t=t(i,c)$ sufficiently large such that
$$
\int_\Omega|\nabla u^t|^2dx=t^2\int_\Omega|\nabla v|^2dx \geq t^2\theta_1c \geq 2\rho_i
$$
and $\widetilde{E}_\tau (u^t)< 0$ uniformly with respect to $\tau$ and $ u\in V_i \cap S_c$.

\end{proof}.

For any $c<\widetilde{c}_i$ and $\tau \in [\frac{1}{2},1]$, define
$$
\nu_{i,\tau,c}=\inf_{\gamma\in\Gamma_{i,c}}\sup_{t\in[0,1];u\in V_i\cap S_c}\widetilde{E}_\tau(\gamma(t,u)),
$$
where
$$
\Gamma_{i,c}:=\{\gamma:[0,1]\times(S_c\cap V_i)\to S_c:\gamma ~~\text{is continuous, old in}~u, 
\gamma(0,u)=u,\gamma(1,u)=\tilde{u}\},
$$
where $\widetilde{u}=u^{t}$ and $t=t(i,c)$ is defined in Lemma \ref{findwideu}.

\begin{remark}\label{remarkforgenus}
  Recall some properties of the cohomological index for spaces with an action of the group $G=\{-1,1\}.$ 
\begin{itemize}
  \item $(i)$ If $G$ acts on $\mathbb{S}^n-1$ via multiplication then $i(\mathbb{S}^n-1)=n.$
  \item $(ii)$ If there exists an equivariant map $X\to Y$ then $i(X)\leq i(Y).$
  \item $(iii)$ Let $X=X_0\cup X_1$ be metrisable and $X_0,X_1\subset X$ be closed $G$-invariant subspaces.
Let $Y$ be a $G$-space and consider a continuous map $\phi:[0,1]\times Y\to X$ such that
each $\phi_t=\phi(t,\cdot):Y\to X$ is equivariant. If $\phi_0(Y)\subset X_0$ and $\phi_1(Y)\subset X_1$ then
$$
i(\mathrm{Im}(\phi)\cap X_0\cap X_1)\geq i(Y).
$$
\end{itemize}
 
Properties $(i)$ and $(ii)$ are standard and hold also for the Krasnoselskii genus. 
Property $(I_3)$ has been proven in~\cite[Corollary 4.11 and Remark 4.12]{bartsch2017natural}. We can now prove Lemma 2.3.
\end{remark}

\begin{lemma}\label{linking}
For any $\gamma \in \Gamma_{i,c}$, there exists $(t,u) \in [0,1] \times (S_c\cap V_i)$ such that 
$$
\gamma(t,u)\in \mathcal{B}_i.
$$
Consequently,  
\begin{align}\label{citauboundfrombelow}
  \nu_{i,\tau, c} >\sup_{u\in V_i \cap S_c} max\{\widetilde{E}_\tau(u), \widetilde{E}_\tau(\widetilde{u})\}, \forall \tau \in [\frac{1}{2},1].
\end{align}

\end{lemma}

\begin{proof}
  
  Let 
  $$
  X=V_{i-1}\times\mathbb{R}^+,\quad X_0=V_{i-1}\times[0,\rho_i],\quad X_1=V_{i-1}\times[\rho_i,+\infty).
  $$
  
  Then $X=X_0\cup X_1.$ Let $P_{i-1}:H_0^1(\Omega)\to V_{i-1}$ be the projection. 
  Define $h_i:S_c\to V_{i-1}\times\mathbb{R}$ by
  $$
  h_i(u)=\left(P_{i-1}(u),\int_\Omega|\nabla u|^2dx\right),
  $$
  and $\phi:[0,1]\times(S_c\cap V_i)\to V_{i-1}\times\mathbb{R}$ by
  $$
  \phi(t,u)=h_i\circ\gamma(t,u),\quad\forall\:(t,u)\in[0,1]\times(S_c\cap V_i).
  $$

 Therefore, $\mathcal{B}_{i}=h_{i}^{-1}(0,\rho_{i}).$ 
Assume by contradiction that
  $$
  \gamma(t,u)\notin\mathcal{B}_i,\quad\forall\:(t,u)\in[0,1]\times(S_c\cap V_i).
  $$
Then
  $$
  \mathrm{Im}(\phi)\cap X_0\cap X_1\subset(V_{i-1}\setminus\{0\})\times\{\rho_i\}.
  $$
  It follows from properties of the genus that
  $$
  \gamma(\mathrm{Im}(\phi)\cap X_0\cap X_1)\leq\gamma((V_{i-1}\setminus\{0\})\times\{\rho_i\})=\dim V_{i-1}.
  $$
  On the other hand, set $\phi_0(u)=\phi(0,u)$ and $\phi_1(u)=\phi(1,u)$ for any $u\in V_i\cap S_c$,
  then
  $$
  \phi_0(S_c\cap V_i)\subset V_{i-1}\times[0,\:\rho_i]=X_0,\quad\phi_1(S_c\cap V_i)\subset V_{i-1}\times[\rho_i,\:+\infty)=X_1.
  $$
  By Property (iii) in Remark~\ref{remarkforgenus},  we have
  $$
  \gamma(\mathrm{Im}(\phi)\cap X_0\cap X_1)\geq\gamma(S_c\cap V_i)=\mathrm{dim}V_i.
  $$

  Hence, we obtain that $dimV_i<dimV_{i-1}$, which contradicts the definition of $V_i$.
 Together with Lemma \ref{lemmaEtauvonBigeqEtauu} and Lemma \ref{findwideu}, we can obtain \eqref{citauboundfrombelow}.
  
\end{proof}

Recall the monotonicity trick applied to $S_c$ as follows, which can be proved in the similar way as Theorem~\ref{Monotonicitytrick}
and \cite[Lemma 3.5]{jeanjean1999existence}. 

\begin{theorem}\label{MonotonicitytrickonSc} (Monotonicity trick).

  Let $I=[\frac{1}{2},1]$. We consider a family $(\widetilde{E}_\tau)_{\tau\in I}$ of $C^1$-functionals on $H_0^1(\Omega)$ of the form
  $$
  \widetilde{E}_\tau(u)=A(u)-\tau B(u),\quad\tau\in I
  $$
  where $B(u)\geq0,\forall u\in H_0^1(\Omega)$ and such that either $A( u) \to \infty$ or $B( u) \to \infty$ as $\| u\| \to \infty.$ 
  We assume that, $\forall \tau \in I$,
  $$
  \nu_{\tau}:=\inf_{\gamma\in\Gamma}\sup_{t\in[0,1]; u_1,u_2 \in S_c}\widetilde{E}_\tau(\gamma(t))
  >\sup_{u_1,u_2 \in S_c}\max\{\widetilde{E}_\tau(u_1),\widetilde{E}_\tau(u_2)\},
  $$
  where 
 $$
 \Gamma=\{\gamma\in C([0,1],S_c),\gamma(0)=u_1,\gamma(1)=u_2\}. 
 $$

  Then, for almost everywhere $\tau\in I$, there is a sequence $\{u_n\}\subset S_c$ such that
  \begin{itemize}
    \item[$(i)$] $\{u_n\}$ is bounded in $H_0^1(\Omega)$;
    \item[$(ii)$] $\widetilde{E}_\tau ( u_n) \to \nu_{\tau}$;
    \item[$(iii)$] $\widetilde{E}_\tau ^{\prime }| _{S_c }( u_n) \to 0$ in $H^{-1}(\Omega)$.
    \end{itemize} 
  
  \end{theorem}

\begin{proposition}\label{multipleaetau1}
  For any $0<c<\widetilde{c}_i$ and almost everywhere $\tau \in [\frac{1}{2},1]$, $\nu_{i,\tau, c}$ is a critical value of $\widetilde{E}_\tau$,
  and there exists a critical $u_\tau$ of $\widetilde{E}_\tau$ constrained on $S_c$ at the level $\nu_{i,\tau, c}$, 
  which solves \eqref{fequationoftauwidetildeE}. 

\end{proposition}

Apply Theorem \ref{MonotonicitytrickonSc} with
$$
A(u)=\frac{1}{2}\int_\Omega|\nabla u|^{2}dx,~~ B(u)=\int_\Omega F(u)dx. 
$$
By Lemma~\ref{linking}, 
we know that for almost everywhere $\tau \in [\frac{1}{2},1]$, there exist a bounded (PS) sequence $\{u_n\} \subset S_c$ 
satisfying $\widetilde{E}_\tau(u_n) \to \nu_{i,\tau, c}$ and $(\widetilde{E}_\tau|_{S_c})'(u_n) \to 0$ as $n \to \infty$.
Hence, up to a subsequence, there exist $u_\tau \in S_c$ such that
\begin{align}\label{multipleuntoutauweak2}
    u_n\rightharpoonup  u_\tau~~\text{ in }H_0^1(\Omega), u_n\to u_\tau~~\text{ in }L^r(\Omega), \forall2~\leq r<2^*,
    u_n\to u_\tau~~\text{a.e. in }\Omega.
\end{align}

We can verify that $u_\tau \in S_c$ is a critical point of $\widetilde{E}_\tau$ constrained on $S_c$.

Let 
$$
\lambda_n=\frac{1}{c}\left(\tau \int_\Omega f(u_n)u_n dx-\int_\Omega|\nabla u_n|^{2}dx\right),
$$
then $\lambda_n$ is bounded and 
\begin{align}\label{Entauprimeto0}
\widetilde{E}_\tau'(u_n)-\lambda_nu_n \to 0~~\text{as}~~ n \to \infty. 
\end{align}

Therefore there exist $\lambda_\tau \in \mathbb R$ such that
\begin{align}\label{multiplelambdanto}
  \lambda_n \to \lambda_\tau
\end{align}

and 
\begin{align}\label{Eprimeto0}
\widetilde{E}_\tau'(u_\tau)+\lambda_\tau u_\tau=0 ~~\text{in}~~H^{-1}(\Omega). 
\end{align}

From~\eqref{multipleuntoutauweak2}, \eqref{Entauprimeto0}, \eqref{multiplelambdanto} and \eqref{Eprimeto0}, we have
$$
u_n \to  u_\tau ~~\text{ in }~~H_0^1(\Omega)
$$
and 
$$
\nu_{i,\tau,c}=\widetilde{E}_\tau(u_\tau).
$$.

\begin{proof}[\textbf{Proof of Theorem~\ref{theorem4A}}]

Similar to Lemma~\ref{matautom1}, we can prove that $\lim_{\tau \to 1^-}\nu_{i,\tau, c}=\nu_{i,1, c}$.
Then we can take $\epsilon$ sufficiently small such that 
$$
 \nu_{i,\tau, c} \leq 2\nu_{i,1, c},~\forall \tau \in [1-\epsilon,1].
$$

Let $\{\tau_n\} \subset [1-\epsilon,1]$ be a sequence  such that $\tau_n \to 1^-$ as $n \to \infty$, and $(\lambda_n, u_n)$ be a
normalized solution of ~\eqref{equationoftau} at the level $\nu_{i,\tau_n, c}$. 
  
On the one hand, we have
\begin{align}\label{fwidetildeEleqnu1}
  \widetilde{E}_{\tau_n}(u_n) = \frac12\int_\Omega|\nabla u_n|^2dx-\tau_n\int_{\Omega}F(u_n)dx
=\nu_{i,\tau_n, c} \leq 2\nu_{i,1, c}. 
\end{align}

On the other hand, we know that any critical point $u$ of $\widetilde{E}_{\tau} |_{S_c}$ satisfies the following Pohozaev identity:
$$
\int_{\Omega}|\nabla u|^2dx-\frac{1}{2}\int_{\partial\Omega}|\nabla u|^2\sigma\cdot nd\sigma=
\tau\left(\frac{N}{2}\int_{\Omega}f(u)udx-N\int_{\Omega}F(u)dx\right). 
$$
  
Note that $\sigma\cdot n>0$ since $\Omega$ is star-shaped with respect to the origin.
Then $u_n$ satisfies
\begin{align}\label{fpohozaevforun}
\int_{\Omega}|\nabla u_n|^2dx>
\tau_n\left(\frac{N}{2}\int_{\Omega}f(u_n)u_ndx-N\int_{\Omega}F(u_n)dx\right). 
\end{align}

Combining~\eqref{fwidetildeEleqnu1} and~\eqref{fpohozaevforun}, we have
$$
(\frac12-\frac{2}{( \alpha -2)N})\int_\Omega|\nabla u_n|^2< (p-2)N \nu_{i,1, c}.
$$

Consequently, we obtain that $u_n$ is bounded in $H_0^1(\Omega)$ uniformly with respect to $n$. 

It follows from ~\eqref{fequationoftauwidetildeE} that
$$
\lambda_n=\frac{1}{c}(\tau_n \int_\Omega f(u_n)u_n dx-\int_\Omega|\nabla u_n|^{2}dx),
$$
and so $\lambda_n$ is bounded.
Therefore, there exists $u_{i,c}$ and $\lambda_{i,c}$ such that
\begin{align}\label{fmultiembeddingandto}
    u_n\rightharpoonup  u_{i,c}~~\text{ in }H_0^1(\Omega), u_n\to u_{i,c}~~\text{ in }L^r(\Omega), \forall~~2~\leq r<2^*,
    u_n\to u_{i,c}~~\text{a.e. in }\Omega.
\end{align}

and
\begin{align}\label{lambdantolambdaic}
\lambda_n \to \lambda_{i,c}.  
\end{align}

Hence, $(\lambda_{i,c}, u_{i,c})$ is a normalized solution of~\eqref{Peps}.
By~\eqref{fequationoftauwidetildeE}, \eqref{fmultiembeddingandto}and  \eqref{lambdantolambdaic}we have
$$
u_n \to u_{i,c}~~\text{in}~~H_0^1(\Omega).
$$
Consequently, $E(u_{i,c})=\nu_{i,1,c}$.

Indeed, by the proof of Lemma~\ref{lemmaEtauvonBigeqEtauu}, Lemma~\ref{linking} and the definition of $\rho_i$, we
can deduce that
$$
\nu_{i,1,c} \geq \frac{1}{4}\rho_i.
$$
Therefore, $\nu_{i,1,c} \to +\infty$ as $i \to \infty$.

Let $K \in \mathbb N$ be such that there are $k$ elements in the set $\{\nu_{i,1,c}|i=1,\cdots, K \}$.
Define 
\begin{align}\label{fdefineck}
c_k:=min\{\widetilde{c}_i|i=1,\cdots,K\}, 
\end{align}

where $\widetilde{c}_i$ is defined by~\eqref{defineci}.
As a result, we obtain that there exist at least $k$ normalized solutions for $0<c<c_k$. 
\end{proof}

\subsection{Multiple normalized solutions when (2) holds}
For any $\tau \in [\frac{1}{2},1]$, we define $\widetilde{E}_\tau(u):S_c \to \mathbb R^N$ by 
$$
\widetilde{E}_\tau(u) = \frac12\int_\Omega|\nabla u|^2dx+\int_\Omega |u|^q-\tau \int_\Omega |u|^p. 
$$

Clearly, the critical point $u$ of $\widetilde{E}_\tau$ on $S_c$ satisfies
\begin{align}\label{equationoftauwidetildeE}
  \begin{cases}-\Delta u+\lambda u=\tau |u|^{p-2}u- |u|^{q-2}u&\text{in }\Omega,\\
    u>0\quad\text{in }\Omega,\quad u=0&\text{on }\partial\Omega.\end{cases}
\end{align}

Define $\varrho_i$ by
\begin{align}\label{defineanotherrhoi}
  \frac{1}{p}C_{N,p}^p \theta_i^{{\frac{N(p-2)-2p}2}}\varrho_i^{\frac{p-2}2}
  +\frac{2^{\frac{-q}{2}}}{q}C_{N,q}^q{\theta_i}^{\frac{N(q-2)-2q}4}\varrho_i^{\frac{q-2}2}
  =\frac{1}{4}  
\end{align}

Since $2<q<p<2^*$, then $\frac{p-2}2, \frac{q-2}2>0$ and $\frac{N(p-2)-2p}4, \frac{N(q-2)-2q}4<0$.
Therefore, $\varrho_i$ is well defined and $\varrho_i \to +\infty$ as $i \to \infty$.

\begin{lemma}\label{lemmaEtauvonBigeqEtauu2}
  Let
  \begin{align}\label{defineci2}
    \widetilde{\alpha}_i=\frac{\varrho_i}{2\theta_i}
  \end{align}
and  $c<\widetilde{\alpha}_i$. Define
$$
\mathcal{B}_i=\{v \in V_{i-1}^\bot \cap S_c:\|\nabla v\|_2^2=\varrho_i \}.
$$ 
Then for any $u\in V_i \cap S_c$, we have
$$
\int_\Omega|\nabla u|^{2}dx<\varrho_i, 
$$
and 
\begin{align}\label{EtauvonBigeqEtauu2}
  \widetilde{E}_\tau(u)< \inf_{v \in \mathcal{B}_i}\widetilde{E}_\tau ~~\text{for any}~~\tau \in[\frac{1}{2},1].
\end{align}

\end{lemma}

\begin{proof}

Since $ u \in V_{i} \cap S_c$, then $\int_\Omega|\nabla u|^2dx\leq \theta_i c<\frac{\varrho_i}{2}$.

For any $u\in V_i \cap S_c$, 

$$
  \begin{aligned}
    \widetilde{E}_\tau(u)&\leq \frac12\int_\Omega|\nabla u|^2dx+\frac{1}{q}\int_\Omega|u|^q\\
    &\leq \frac{1}{2}\theta_i c+\frac{1}{q}C_{N,q}^q c^{\frac{2q-N(q-2)}4}(\theta_i c)^{\frac{N(q-2)}4}\\
    &< (\frac{1}{4}+\frac{2^{\frac{-q}{2}}}{q}C_{N,q}^q{\theta_i}^{\frac{N(q-2)-2q}4}\varrho_i^{\frac{q-2}2})\varrho_i.
    \end{aligned}
$$

On the other hand, for any $v\in \mathcal{B}_i$, 
$$
  \begin{aligned}
    \widetilde{E}_\tau(v)&\geq \frac12\int_\Omega|\nabla v|^2dx-\frac{1}{p}\int_\Omega|v|^p\\
    &\geq \frac{1}{2}\varrho_i-\frac{1}{p}C_{N,p}^p c^{\frac{2p-N(p-2)}4}\varrho_i^{\frac{N(p-2)}4}\\
    &= \left(\frac{1}{2}-\frac{1}{p}C_{N,p}^p \theta_i^{{\frac{N(p-2)-2p}4}}\varrho_i^{\frac{p-2}2}\right)\varrho_i.
    \end{aligned}
$$

It follows from the definition of $\varrho_i$ that $\widetilde{E}_\tau(u)< \inf_{v \in \mathcal{B}_i}\widetilde{E}_\tau$.

\end{proof}

Direct calculations yield that
$$
  \begin{aligned}
\widetilde{E}_\tau(u^t) &=\frac12\int_\Omega|\nabla u^t|^2dx+\int_\Omega |u^t|^q-\tau \int_\Omega |u^t|^p \\
&\leq \frac{t^2}{2}\int_\Omega|\nabla u|^2dx+t^{N(\frac{q-2}{2})}\int_\Omega| u|^{q}dx-\frac{t^{N(\frac{p-2}{2})}}{2}\int_\Omega| u|^{p}dx\\
& \leq \frac{t^2}{2} \theta_i c +t^{N(\frac{q-2}{2})}C_{N,q}^q c^{\frac{2q-N(q-2)}4}(\theta_i c)^{\frac{N(q-2)}4}
-\frac{t^{N(\frac{p-2}{2})}}{2} c^{\frac{p}{2}}|\Omega|^{\frac{2-p}{2}}. 
    \end{aligned}
$$

Since $N(\frac{p-2}{2})>\max\{N(\frac{q-2}{2}), 2\}$, then $\widetilde{E}_\tau(u^t) \to -\infty$ as $t \in +\infty$. 
Similar to Lemma \ref{findwideu}, we can prove there exists $t=t(i,c)$ such that for $\widetilde{u}=u^t$, 
$\int_\Omega|\widetilde{u}|^2=\int_\Omega|u|^2=c$, 
$\int_\Omega|\nabla \widetilde{u}|^2dx \geq 2\rho_i$ and
$\widetilde{E}_\tau (\widetilde{u})< 0$ uniformly with respect to $\tau$ and $ u\in V_i \cap S_c$.

For any $c<\widetilde{\alpha}_i$ and $\tau \in [\frac{1}{2},1]$, define
$$
\nu_{i,\tau,c}=\inf_{\gamma\in\Gamma_{i,c}}\sup_{t\in[0,1];u\in V_i\cap S_c}\widetilde{E}_\tau(\gamma(t,u)),
$$
where
$$
\Gamma_{i,c}:=\{\gamma:[0,1]\times(S_c\cap V_i)\to S_c:\gamma ~~\text{is continuous, old in}~u, 
\gamma(0,u)=u,\gamma(1,u)=\tilde{u}\}. 
$$

\begin{proposition}\label{multipleaetau}
  For any $0<c<\widetilde{\alpha}_i$ and almost everywhere $\tau \in [\frac{1}{2},1]$, $\nu_{i,\tau, c}$ is a critical value of $\widetilde{E}_\tau$,
  and there exists a critical $u_\tau$ of $\widetilde{E}_\tau$ constrained on $S_c$ at the level $\nu_{i,\tau, c}$, which solves \eqref{equationoftauwidetildeE}.

\end{proposition}

\begin{proof}
    The proof is similar to the proof of Proposition~\ref{multipleaetau1}. 
\end{proof}

\begin{proof}[\textbf{Proof of Theorem~\ref{theorem4B}}]

Similar to Lemma~\ref{matautom1}, we can prove that $\lim_{\tau \to 1^-}\nu_{i,\tau, c}=\nu_{i,1, c}$.
Then we can take $\epsilon$ sufficiently small such that 
$$
 \nu_{i,\tau, c} \leq 2\nu_{i,1, c},~\forall \tau \in [1-\epsilon,1].
$$

Let $\{\tau_n\} \subset [1-\epsilon,1]$ be a sequence  such that $\tau_n \to 1^-$ as $n \to \infty$, and~\eqref{equationoftau} a
normalized solution $(\lambda_n, u_n)$ with energy $\nu_{i,\tau_n, c}$. 
  
On the one hand, 
\begin{align}\label{widetildeEleqnu1}
  \widetilde{E}_{\tau_n}(u_n) = \frac12\int_\Omega|\nabla u_n|^2dx+\frac{1}{q}\int_\Omega|u_n|^q-\frac{\tau_n}{p}\int_\Omega|u_n|^p
=\nu_{i,\tau_n, c} \leq 2\nu_{i,1, c}. 
\end{align}

On the other hand, we know that any critical point $u$ of $\widetilde{E}_\tau|_{S_c}$ satisfies the following Pohozaev identity:
$$
\int_{\Omega}|\nabla u|^2dx-\frac{1}{2}\int_{\partial\Omega}|\nabla u|^2\sigma\cdot nd\sigma=
\tau(\frac{1}{2}-\frac{1}{p})N\int_\Omega |u|^p-(\frac{1}{2}-\frac{1}{q})N\int_\Omega |u|^q. 
$$
  
Note that $\sigma\cdot n>0$ since $\Omega$ is star-shaped with respect to the origin.
Then $u_n$ satisfies
\begin{align}\label{pohozaevforun}
\int_{\Omega}|\nabla u_n|^2dx+(\frac{1}{2}-\frac{1}{q})N\int_\Omega |u_n|^q>
\tau_n (\frac{1}{2}-\frac{1}{p})N\int_\Omega |u_n|^p. 
\end{align}

Combining~\eqref{widetildeEleqnu1} and~\eqref{pohozaevforun}, we have
$$
\frac{N(p-2)-4}{4}\int_\Omega|\nabla u_n|^2dx+\frac{N(p-q)}{2q}\int_\Omega |u_n|^q < (p-2)N \nu_{i,1, c}.
$$

Consequently, we obtain that $u_n$ is bounded in $H_0^1(\Omega)$ uniformly with respect to $n$. 

It follows from ~\eqref{equationoftauwidetildeE} that
$$
\lambda_n=\frac{1}{c}(\frac{\tau_n}{p}\int_\Omega |u_n|^pdx-\frac{1}{q}\int_\Omega |u_n|^qdx-\int_\Omega|\nabla u_n|^{2}dx),
$$
and so $\lambda_n$ is bounded.
Therefore, there exists $u_{i,c}$ and $\lambda_{i,c}$ such that
\begin{align}\label{multiembeddingandto}
    u_n\rightharpoonup   u_{i,c}~~\text{ in }H_0^1(\Omega), u_n\to  u_{i,c}~~\text{ in }L^r(\Omega), \forall2~\leq r<2^*,
    u_n\to  u_{i,c}~~\text{a.e. in }\Omega,
\end{align}

and
\begin{align}\label{untouic}
\lambda_n \to \lambda_{i,c}.
\end{align}

Hence, $(\lambda_{i,c}, u_{i,c})$ is a normalized solution of~\eqref{Peps}.
By~\eqref{equationoftauwidetildeE},\eqref{untouic} and~~\eqref{multiembeddingandto}, we have
$$
    u_n \to u_{i,c}~~\text{in}~~H_0^1(\Omega). 
$$

Consequently, $E(u_{i,c})=\nu_{i,1,c}$.

Indeed, by the proof of Lemma~\ref{lemmaEtauvonBigeqEtauu}, Lemma~\ref{linking} and the definition of $\rho_i$, we
can deduce that
$$
\nu_{i,1,c} \geq \frac{1}{4}\rho_i.
$$
Therefore, $\nu_{i,1,c} \to +\infty$ as $i \to \infty$.

Let $K \in \mathbb N$ be such that there are $k$ elements in the set $\{\nu_{i,1,c}|i=1,\cdots, K \}$.
Define 
\begin{align}\label{up-uqdefineck}
\alpha_k:=min\{\widetilde{\alpha}_i|i=1,\cdots,K\}, 
\end{align}

where $\widetilde{\alpha}_i$ is defined by~\eqref{defineci2}.
As a result, we obtain that there exist at least $k$ normalized solutions for $0<c<\alpha_k$. 
\end{proof}

\section{Non-radial sign-changing normalized solutions on a ball}

In this section, we consider the existence of non-radial sign-changing normalized solutions of~\eqref{Peps} when $f$ satisfies $(f_1)$
, $\Omega=B$ is a ball and $N \geq 4$.

Let $2\leq \kappa \leq \frac{N}{2}$ be a fixed integer different from $\frac{N-1}{2}$.
The action of
$$
G=O(\kappa) \times  O(\kappa) \times  O(N-2\kappa)
$$
on $H_0^1(\Omega)$ is defined by
$$
gu(x):=u(g^{-1}x).
$$
Let $\eta$ be the involution defined on $\mathbb R^N=\mathbb R^m \oplus \mathbb R^m \oplus \mathbb R^{N-2m}$ by
$$
\eta(x_1, x_2, x_3):=(x_2, x_1, x_3).
$$
The action of $H={id, \eta}$ on $H_{0,G}^1(B)$ is defined by 

\begin{align}
  hu(x):=
\begin{cases}u(x), \qquad \qquad h=id,
    \\ -u(h^{-1}x), \qquad  h=\eta.\end{cases}
  \end{align}

Set 
$$
K=\{u \in H_{0,G}^1(B): ~hu=u, \forall h \in H\}.
$$
Note that $B$ is a radial domain, thus $K$ is well-defined. 
It is not difficult to verify that $0$ is the only radial function in $K$, and the embedding $K \hookrightarrow L^p(\Omega)$ is compact. 

Define $\rho_i$ by
\begin{align}\label{definevarrhoi}
  \zeta (C_{N,\alpha}^\alpha \theta_i^{{\frac{N(\alpha-2)-2\alpha}4}}\rho_i^{\frac{\alpha-2}2}
  +C_{N,\beta}^\beta \theta_i^{{\frac{N(\beta-2)-2\beta}4}}\rho_i^{\frac{\beta-2}2})=\frac{1}{2}.
\end{align}

Similar to the definition of ~\eqref{definerhoi}, we know that $\rho_i$ is well defined and $\rho_i \to +\infty$ as $i \to \infty$.

\begin{lemma}\label{lemmaEtauvonBigeqEtauunonradial}
  Let
  \begin{align}\label{definecinonradial}
    \widetilde{c}_i=\frac{\rho_i}{2\theta_i}
  \end{align}
and  $c<\widetilde{c}_i$. where $\rho_i$ is defined in ~\eqref{definerhoi}. Define
$$
\widetilde{\mathcal{B}}_i=\{v \in V_{i-1}^\bot \cap S_c \cap K:\|\nabla v\|_2^2=\rho_i \}.
$$ 
Then for any $u\in V_i \cap S_c$, we have
$$
\int_\Omega|\nabla u|^{2}dx<\rho_i, 
$$
and 
\begin{align}\label{EtauvonBigeqEtauunonradial}
  \widetilde{E}_\tau(u) < \inf_{v \in \mathcal{B}_i}\widetilde{E}_\tau ~~\text{for any}~~\tau \in[\frac{1}{2},1].
\end{align}

\end{lemma}

\begin{proof}
  
  Since $ u \in V_{i} \cap S_c \cap K$, then $\int_\Omega|\nabla u|^2dx\leq \theta_i c<\frac{\rho_i}{4}$.

  For any $u\in V_i \cap S_c \cap K$, 
  $$
  \widetilde{E}_\tau(u) \leq \frac12\int_\Omega|\nabla u|^2dx
  \leq \frac{1}{2}\theta_i c <\frac{\rho_i}{4}.
  $$
 
On the other hand, for any $v \in \widetilde{\mathcal{B}}_i$,

$$
  \begin{aligned}
    \widetilde{E}_\tau(v) &\geq \frac12\int_\Omega|\nabla v|^2dx-\frac{1}{2}\zeta (\int_\Omega|v|^\alpha+\int_\Omega|v|^\beta)\\
  &\geq \frac{1}{2}\rho_i-\frac{1}{2}\zeta (C_{N,\alpha}^\alpha c^{\frac{2\alpha-N(p-2)}4}\rho_i^{\frac{N(\alpha-2)}4}
  +C_{N,\beta}^\beta c^{\frac{2\beta-N(p-2)}4}\rho_i^{\frac{N(\beta-2)}4} )\\
  &= \left(\frac{1}{2}-\frac{1}{2}\zeta (C_{N,\alpha}^\alpha \theta_i^{{\frac{N(\alpha-2)-2\alpha}4}}\rho_i^{\frac{\alpha-2}2}
  +C_{N,\beta}^\beta \theta_i^{{\frac{N(\beta-2)-2\beta}4}}\rho_i^{\frac{\beta-2}2})\right)\rho_i
  \end{aligned}
$$

It follows from the definition of $\rho_i$ that $\widetilde{E}_\tau(u)< \inf_{v \in \widetilde{\mathcal{B}}_i}\widetilde{E}_\tau$.
  
  \end{proof}

  For any $c<\widetilde{c}_i$ and $\tau \in [\frac{1}{2},1]$, define
  $$
  \widetilde{\nu}_{i,\tau,c}=\inf_{\gamma\in\Gamma_{i,c}}\sup_{t\in[0,1];u\in V_i\cap S_c \cap K}\widetilde{E}_\tau(\gamma(t,u)),
  $$
  where
  $$
  \widetilde{\Gamma}_{i,c}:=\{\gamma:[0,1]\times(S_c\cap V_i \cap K)\to S_c:\gamma ~~\text{is continuous, old in}~u, 
  \gamma(0,u)=u,\gamma(1,u)=\tilde{u}\},
  $$
and $\tilde{u}$ is same as that defined in Section 4.

  \begin{proof}[\textbf{Proof of Theorem~\ref{theorem5}}]

We can prove the verison of the monotonicity trick applied to $S_c \cap K$ in the same way as the proof of Theorem \ref{Monotonicitytrick}. 
Furthermore, similar to the proof of Lemma~\ref{linking}, Proposition~\ref{multipleaetau} and Theorem\ref{theorem4A}, 
 we can prove that 
  \begin{align}\label{widetildecitauboundfrombelow}
    \widetilde{\nu}_{i,\tau, c} >
     \sup_{u \in V_i \cap S_c \cap K} max\{\widetilde{E}_\tau(u), \widetilde{E}_\tau(\widetilde{u})\}, \forall \tau \in [\frac{1}{2},1],
  \end{align}
and obtain the critical point in the same way. 

\end{proof}

\medskip

\textbf{Acknowledgements}
I would like to thank C. Li and S. J. Li for fruitful disscussions and constant support during the development of this work.

\end{document}